# MODERATE DEVIATION PROBABILITIES FOR OPEN CONVEX SETS: NONLOGARITHMIC BEHAVIOR

By Uwe Einmahl[1] and James Kuelbs[2]

*Vrije Universiteit Brussel and University of Wisconsin*

Precise asymptotics for moderate deviation probabilities are established for open convex sets in both the finite- and infinite-dimensional settings. Our results are based on the existence of dominating points for these sets, a related representation formula, and asymptotics for the integral term in this formula.

**1. Introduction.** Let $X, X_1, X_2, \ldots$ be independent, identically distributed random vectors where $\mathcal{L}(X) = \mu$, and $\mu$ is a Borel probability measure on the real separable Banach space $B$. Let $S_n = \sum_{j=1}^n X_j$ and assume $\mathcal{L}(S_n/n^{1/2})$ converges weakly. Then the limit law $\gamma$ is necessarily Gaussian with mean zero, and $\mu$ also has mean zero. Let $\{b_n\}$ be a positive sequence such that

$$b_n/n^{1/2} \to \infty \quad \text{and} \quad b_n/n \to 0. \tag{1.1}$$

Here we study the asymptotic behavior of $\{P(S_n/b_n \in A)\}$ under (1.1). These probabilities are frequently called moderate deviation probabilities, and there is a long history of such results in the finite-dimensional setting. There are also results in the infinite-dimensional setting, but only at the logarithmic level. In particular, the results by Borovkov and Mogul'skii [6] and by de Acosta [9] are of this type.

Let $B^*$ denote the topological dual space of $B$ and define

$$\hat{\mu}(f) = \int_B e^{f(x)} \, d\mu(x), \qquad f \in B^*,$$

$$\hat{\gamma}(f) = \int_B e^{f(x)} \, d\gamma(x), \qquad f \in B^*. \tag{1.2}$$

Received August 2002; revised April 2003.
[1] Supported in part by an FWO-grant.
[2] Supported in part by NSF Grant DMS 00-71-700.

Since $\gamma$ is centered Gaussian, $\hat{\gamma}(f) = \exp\{\sigma_f^2/2\}$, where $\sigma_f^2 = \int_B f^2(x)\,d\gamma(x)$. Furthermore, it is well known that the rate function

(1.3) $$\lambda_\gamma(x) = \sup_{f \in B^*} [f(x) - \log \hat{\gamma}(f)], \qquad x \in B,$$

is given by

(1.4) $$\lambda_\gamma(x) = \begin{cases} \|x\|_\gamma^2/2, & \text{if } x \in H_\gamma \subset B, \\ +\infty, & \text{otherwise.} \end{cases}$$

Here $H_\gamma$ is the Hilbert space generating $\gamma$ on $B$, that is, the completion of $S(B^*)$ where $S: B^* \to B$ is given by the integral

(1.5) $$Sf = \int_B x f(x)\,d\gamma(x), \qquad f \in B^*,$$

in the norm determined by the inner product $\langle Sf, Sg \rangle_\gamma = \int_B f(x)g(x)\,d\gamma(x)$. Since $\gamma$ has moments of all order, $Sf$ exists as a Bochner integral. Further details can be found in Lemma 2.1 of [14].

Let $\overline{D}$ denote the $B$-closure of $D$ and $\partial D$ the boundary of $D$. Throughout we assume that

(1.6) 
    (i)   $D$ is an open convex subset of $B$.
    (ii)  $D \cap H_\gamma \neq \phi$.
    (iii) $0 \notin \overline{D}$.

Since $\int_B e^{t\|x\|}\,d\gamma(x)$ for all $t > 0$, then [12], Theorem 1, implies $D$ has a unique dominating point with respect to $\gamma$ (see also [15] and [16]). That is, there exists a unique point $a_0 \in \partial D$ such that

(1.7)
    (i)   $\lambda_\gamma(a_0) = \inf_{x \in D} \lambda_\gamma(x) = \inf_{x \in \overline{D}} \lambda_\gamma(x) < \infty$.
    (ii)  For some $g \in B^*$ we have $D \subset \{x : g(x) \geq g(a_0)\}$.
    (iii) $\lambda_\gamma(a_0) = g(a_0) - \log \hat{\gamma}(g)$ and
    (iv) $a_0 = \int_B x \exp\{g(x) - \log \hat{\gamma}(g)\}\,d\gamma(x)$, where the integral exists as a Bochner integral.

Furthermore, if we apply the Hahn–Banach theorem and take $f \in B^*$ such that

(1.8) $$\sup_{\{z : \lambda_\gamma(z) \leq \lambda_\gamma(a_0)\}} f(z) = f(a_0) < f(x) \qquad \forall x \in D,$$

then [12], Theorem 1, implies there exists a unique $t_0 > 0$ such that $g = t_0 f$, satisfies (1.7)(ii)–(iv).

In [6], Borovkov and Mogul'skii prove the following result.

THEOREM A. *Let $X, X_1, X_2, \ldots$ be i.i.d. $B$-valued with $\mathcal{L}(S_n/n^{1/2})$ converging weakly to the Gaussian measure $\gamma$ and assume $D$ is an open convex subset of $B$. If $\{b_n\}$ satisfies (1.1) and*

(1.9) $$E(e^{t|f(X)|}) < \infty, \qquad 0 < |t| < t_f, f \in B^*,$$

*then*

$$\lim_{n \to \infty} nb_n^{-2} \log P(S_n/b_n \in D) = - \inf_{x \in D} \lambda_\gamma(x), \tag{1.10}$$

where $\lambda_\gamma$ is given by (1.3).

Under additional integrability assumptions, a full moderate deviation principle for open and closed sets (in the sense of Varadhan) is established for $\{\mathcal{L}(S_n/b_n)\}$ by de Acosta in [9]. In addition, the papers [8] and [17] deal with necessary and sufficient conditions for the upper bound for closed sets in the large deviation principle for various sequences $\{b_n\}$. These results are at the logarithmic level and are quite different from what we establish in the results that follow.

Our interest here is to seek refinements of Theorem A which allow us to study the behavior of $P(S_n/b_n \in D)$ directly, not merely at the logarithmic level. This will be done via a representation formula, which is elementary to establish once one has dominating points, and is the analogue of a similar formula in the large deviation setting. This representation formula becomes useful for moderate deviation probabilities when, in addition to $b_n/n^{1/2} \to \infty$, we also assume $b_n/n^{2/3} \to 0$. What we find is that is in this range, the moderate deviation probabilities are much the same as those when $\mathcal{L}(X) = \gamma$. This is standard in $\mathbb{R}$, but less well understood in the vector space setting.

Our results depend on the shape of $D$ at the dominating point $a_0 \in \partial D$, and the difficult part of our arguments involves establishing the appropriate lower bounds. For upper bounds, replacing $D$ by a half-space is frequently good enough provided $D$ is sufficiently round at $a_0$.

As usual, $a_n \sim b_n$ means $\lim_n a_n/b_n = 1$.

THEOREM 1. *Let $X, X_1, X_2, \ldots$, be i.i.d. $B$-valued random vectors, where $B$ is a separable Banach space, and set $S_n = \sum_{j=1}^n X_j$. Assume $\{S_n/n^{1/2}\}$ converges weakly to a nondegenerate probability measure $\gamma$ on $B$, and that $\{b_n\}$ is a sequence of positive constants such that*

$$b_n/n^{1/2} \to \infty \quad \text{and} \quad b_n/n^{2/3} \to 0. \tag{1.11}$$

*In addition, assume that $D$ satisfies (1.6), (1.9) holds, $a_0$ is the unique dominating point for $(D, \gamma)$, and $g = t_0 f$ is as in (1.7) and (1.8). Then*

$$\begin{aligned} P(S_n/b_n \in D) \\ \sim \exp\{-n^{-1}b_n^2 \lambda_\gamma(a_0)\} E[\exp\{-g(T_n - E(T_n))\} I\{T_n \in b_n^2 D/n\}], \end{aligned} \tag{1.12}$$

*where $T_n = \frac{b_n}{n} \sum_{j=1}^n Z_{n,j}$, and $Z^{(n)}, Z_{n,1}, Z_{n,2}, \ldots, Z_{n,n}$ are i.i.d. with $Z^{(n)}$ being a $B$-valued random variable such that*

$$\frac{d\mathcal{L}(Z^{(n)})}{d\mu}(x) = \exp\{g(b_n x/n)\}/\hat{\mu}(b_n g/n), \tag{1.13}$$

and

(1.14) $$E(Z^{(n)}) = (b_n/n)a_0 + O(b_n^2/n^2).$$

Furthermore,

(1.15) $$\limsup_{n\to\infty} n^{-1/2} b_n P(S_n/b_n \in D) \exp\{n^{-1} b_n^2 \lambda_\gamma(a_0)\} \leq (2\pi\sigma_g^2)^{-1/2},$$

where $\sigma_g^2 = E(g^2(X))$.

To establish lower bounds comparable to (1.15) we need the following definition.

DEFINITION 1. Assume (1.6) and let $a_0$ be the unique dominating point of $D$ with respect to $\gamma$. Then, $D$ contains slices whose diameters near $a_0$ dominate the function $\tau(s)$ if for some $f \in B^*$ satisfying (1.8) there exists $x_0 \in B$, and $\delta > 0$ such that $f(x_0) > 0$, and

(1.16) $$\{y + sx_0 : f(y) = 0, \|y\| \leq \tau(s), 0 < s \leq \delta\} \subset D - a_0.$$

Our first lower bound result is the following theorem.

THEOREM 2. Let $\{b_n\}$ satisfy (1.11) and assume $X, X_1, X_2, \ldots,$ and $\{S_n\}$ satisfy the assumptions of Theorem 1. Also assume

(1.17) $$E(\|X\|^3 e^{t|f(X)|}) < \infty, \qquad 0 < |t| < t_f, f \in B^*$$

and that $D$ satisfies (1.6). Let $a_0$ be the unique dominating point for $(D, \gamma)$ and $g = t_0 f$ be as in (1.7) and (1.8). If $\{\sum_{j=1}^n (Z_{n,j} - E(Z_{n,j}))/n^{1/2}\}$ is bounded in probability, where $Z_{n,1}, Z_{n,2}, \ldots, Z_{n,n}$ are as in Theorem 1 and $D$ contains slices whose diameters near $a_0$ dominate the function $\tau(s) = \beta(s|\log s|)^{1/2}$, $\beta > 0$, then

(1.18) $$\liminf_{n\to\infty} n^{-1/2} b_n P(S_n/b_n \in D) \exp\{n^{-1} b_n^2 \lambda_\gamma(a_0)\} > 0.$$

We note that if $B$ is a Hilbert space or more generally a type 2 Banach space, then the condition on stochastic boundedness follows easily from (1.17). Moreover in the Hilbert space case, Theorem 2 can be improved as follows.

THEOREM 3. Let $\{b_n\}$ satisfy (1.11) and assume that $X, X_1, X_2, \ldots,$ are i.i.d. random vectors taking values in a separable Hilbert space $H$ with (1.17) holding and $E(X) = 0$. Let $D$ satisfy (1.6) and assume $a_0$ is the unique dominating point for $(D, \gamma)$. If $D$ contains slices whose diameter near $a_0$ dominate the power function $\tau(s) = \beta s^{1/2}, \beta > 0$, then (1.18) holds.

We note that Theorem 3, in particular, applies if $D$ is a ball in a Hilbert space satisfying (1.6). (This follows, for instance, from the proof of Theorem 3, [16].)

If $H$ is $\mathbb{R}^d$, then we can obtain more precise estimates of these moderate deviation probabilities. This is our next result.

THEOREM 4. *Let $\{b_n\}$ satisfy* (1.11) *and assume $X, X_1, X_2, \ldots,$ are i.i.d. $\mathbb{R}^d$-valued with $\mathcal{L}(S_n/n^{1/2})$ converging weakly to a Gaussian measure $\gamma$ on $\mathbb{R}^d$ with the support of $\gamma$ all of $\mathbb{R}^d$. Also assume*

$$(1.19) \qquad E(e^{t|\langle f, X\rangle|}) < \infty, \qquad 0 < t < t_f, f \in \mathbb{R}^d,$$

*and let $D$ be as in Theorem* 3. *Then*

$$(1.20) \qquad \lim_{n\to\infty} P(S_n \in b_n D)/P(G \in n^{-1/2} b_n D) = 1,$$

*where $\mathcal{L}(G) = \gamma$.*

If $D$ is a ball we can extend the last result to infinite-dimensional Hilbert space valued random vectors.

THEOREM 5. *Let $H$ be a separable Hilbert space and let $X, X_1, X_2, \ldots,$ be i.i.d. $H$-valued random vectors as in Theorem* 3. *Let $G$ be a Gaussian random vector on $H$ with $\mathcal{L}(G) = \gamma$. If $D = \{x : \|x - a\| < R\}$ is a ball in $H$ satisfying* (1.6.ii) *and* (1.6.iii), *where $\|\cdot\|$ is the Hilbert space norm on $H$, then we have, for any sequence $\{b_n\}$ satisfying* (1.11),

$$(1.21) \qquad \lim_{n\to\infty} P(S_n \in b_n D)/P(G \in n^{-1/2} b_n D) = 1.$$

*Furthermore, both probabilities are asymptotically equivalent to the quantity*

$$(1.22) \quad (2\pi \sigma_g^2 b_n^2/n)^{-1/2} \exp\{-n^{-1} b_n^2 \lambda_\gamma(a_0)\} \int_0^\infty e^{-s} P(\|G_2\|^2 \leq 2sbR^2)\,ds,$$

*where $a_0$ is the unique dominating point for $(D, \gamma)$ and $g = t_0 f$ is as in* (1.7) *and* (1.8), $1/b = g(a - a_0)$, $\sigma_g^2 = E(g^2(X))$, *and $G_2 = G - G_1$ is a centered Gaussian random vector on $H$ with $G_1 = g(G)E(Gg(G))/\sigma_g^2$.*

The remaining part of the paper is organized as follows: We prove Theorem 1 in Section 2. Then we prove Theorem 2 in Section 3, where we use modifications of arguments from [12] when $b_n = n$. The proof of Theorem 3 appears in Sections 4 and 5, and follows from Proposition 1, which depends on a Berry–Esseen result for $U$-statistics from [1]. When $b_n = n$, the analogue of Proposition 1 in [12] was proved via a Berry–Esseen result for $U$-statistics due to van Zwet [19], but this result is no longer applicable when $\lim_n b_n/n = 0$. Hence, we developed a direct approach (independent of $U$-statistics) for proving Proposition 1 in this setting. A refinement of this method allowed us also to eventually prove Proposition 2, which is crucial for obtaining the precise results for balls in Hilbert space given in Theorem 5. Subsequent discussions with V. Bentkus made us aware of some recent improvements of van Zwet's Berry–Esseen inequality for $U$-statistics which appear in [1] and [2]. Once we had these results at our disposal, the proof of Proposition 1 now follows along lines similar to the companion result in [12]. However, the exact asymptotics given in Proposition 2 do not follow in this manner and our "direct" method is still needed for obtaining Theorem 5. Theorem 4 is proved in Section 6, and Theorem 5 in Sections 7 and 8. Both of these theorems provide exact asymptotics for certain open convex sets. In view of relation (1.12) this requires a precise comparison of

$$E[\exp\{-g(T_n - E(T_n))\}I\{T_n \in b_n^2 D/n\}]$$

with a corresponding expectation involving Gaussian random vectors.

To that end, we use, in the finite-dimensional case, an estimate of the convergence speed in the multivariate central limit theorem due to Zaitsev [20] among other tools.

The proof of Theorem 5 (open balls in Hilbert space) is based on Proposition 2 in Section 7. One can rewrite the above expectation as an integral with respect to the two-dimensional distribution of $(\|S_n/n^{1/2}\|^2, f(S_n/n^{1/2}))$, where $f: H \to \mathbb{R}$ is a continuous linear functional. We then show that this distribution is close to that of $(\|Y_n\|^2, f(Y_n))$, where $Y_n$ is an appropriate Gaussian random vector. To accomplish this we need, among other things, a local limit result for a smoothed and truncated version of $(\|S_n/n^{1/2}\|^2, f(S_n/n^{1/2}))$, see Lemma 18. To prove this result we use an adaptation of the characteristic function method for proving Berry–Esseen type results in Hilbert space. For a nice account of this method refer to [3].

**2. Proof of Theorem 1.** The proof of Theorem 1 proceeds with a sequence of lemmas. Throughout this section the conditions of Theorem 1 are assumed. Also note that since $D$ satisfies (1.6), and $g = t_0 f$ relates to $a_0$ as in (1.7) and (1.8), we have $\sigma_g^2 = E(g^2(X)) > 0$.

LEMMA 1. *Let $Z^{(n)}$ be defined as in (1.13), where $\mu = \mathcal{L}(X)$ and $g = t_0 f \in B^*$ is related to the dominating point $a_0$ is in (1.7) and (1.8). Then*

$$(2.1) \qquad E(Z^{(n)}) = \frac{b_n}{n} a_0 + O\left(\frac{b_n^2}{n^2}\right),$$

$$(2.2) \qquad \sigma_{g,n}^2 := E(g^2(Z^{(n)} - E(Z^{(n)}))) = \sigma_g^2 + O\left(\frac{b_n}{n}\right) \quad \text{and}$$

$$(2.3) \qquad a_0 = E(Xg(X)).$$

PROOF. First observe that since $\{S_n/n^{1/2}\}$ converges weakly to $\gamma$, then $\mu$ and $\gamma$ must have the same covariance function, $\gamma$ is a mean zero Gaussian measure, $E\|X\|^{2-\varepsilon} < \infty$ for all $\varepsilon > 0$, and $E(X) = 0$. Hence,

$$(2.4) \qquad E(Xg(X)) = \int_B xg(x)\,d\gamma(x).$$

If $h = Sg$, $S$ given by (1.5), and $\mathcal{L}(Y) = \gamma$, then the Cameron–Martin formula implies

$$(2.5) \qquad h = E(Y+h) = \int_B (x+h)\,d\gamma(x) = \int_B x e^{g(x) - \sigma_g^2/2}\,d\gamma(x).$$

Hence, (2.4), (2.5) and (1.7)(iv) imply $h = Sg = a_0$ and (2.3) holds.

To verify (2.1) we first observe that since $E(X) = 0$,

$$(2.6) \quad \begin{aligned} \hat{\mu}\left(\frac{b_n g}{n}\right) &= E(e^{g(b_n X/n)}) \\ &= E\left(1 + \frac{b_n}{n}g(X) + \frac{1}{2}\left(\frac{b_n}{n}\right)^2 g^2(X) e^{\theta(g(b_n X/n))}\right), \\ &= 1 + \frac{b_n^2}{2n^2} E(g^2(X) e^{\theta b_n g(X)/n}), \end{aligned}$$

where $|\theta| \leq 1$ by Taylor's formula. Since $E(\|X\|^{2-\varepsilon}) < \infty$ and (1.9) is assumed, Hölder's inequality implies $E(\|X\| e^{|g(b_n X/n)|})$ exists for $n$ sufficiently large. Thus, $E(X e^{g(b_n X/n)})$ exists as a Bochner integral for such $n$ and since $E(X) = 0$, we have that

$$(2.7) \quad \begin{aligned} \left\| E(X e^{g(b_n X/n)}) - \frac{b_n}{n}a_0 \right\| &= \left\| E(X e^{g(b_n X/n)}) - E\left(X\left(1 + \frac{b_n g(X)}{n}\right)\right) \right\| \\ &\leq \frac{b_n^2}{2n^2} E(\|X\| g^2(X) e^{|g(b_n X/n|}). \end{aligned}$$

In (2.7) we used (2.3), and if $n$ is large enough, the integral $E(Xg^2(X) e^{|g(b_n X/n|})$ exists as a Bochner integral by an argument similar to that mentioned prior to (2.7). Since $E(Z^{(n)}) = E(X e^{g(b_n X/n)})/\hat{\mu}(b_n g/n)$, we have (2.1) because $b_n/n \to 0$ and the dominated convergence theorem applies. To prove (2.2), we observe

$$(2.8) \quad \begin{aligned} E(g^2(Z^{(n)})) &= E(g^2(X) e^{g(b_n X/n)})/\hat{\mu}(b_n g/n) \\ &= E\left(g^2(X) + \frac{b_n}{n}g^3(X) e^{\theta g(b_n X/n)}\right) \Big/ \hat{\mu}\left(\frac{b_n g}{n}\right), \end{aligned}$$

where $|\theta| \leq 1$. Hence, by (2.6), (2.1) and the dominated convergence theorem, $b_n/n \to 0$ implies (2.2). $\square$

LEMMA 2. *If* (1.1) *holds, then*

$$\lim_{n \to \infty} nb_n^{-2} \log E(e^{f(b_n S_n/n)}) = E(f^2(X))/2 \quad (2.9)$$

*for all* $f \in B^*$. *Furthermore,*

$$b_n^2 n^{-1}[E(f^2(X))/2 - nb_n^{-2} \log E(e^{f(b_n S_n/n)})] = O(b_n^3/n^2). \quad (2.10)$$

PROOF. Since $X, X_1, X_2, \ldots$ are i.i.d., the argument for (2.6) implies

$$\log E(e^{f(b_n S_n/n)}) = n \log E(e^{f(b_n X/n)})$$
$$= n \log\left(1 + \frac{b_n^2}{2n^2}E(f^2(X)) + \frac{b_n^3}{6n^3}E(f^2(X)e^{\theta f(b_n X/n)})\right),$$

where $|\theta| \leq 1$ and $f \in B^*$. Hence, by the dominated convergence theorem and that $\log(1+x) = x + O(x^2)$ as $x \to 0$, we see

$$\log E(e^{f(b_n S_n/n)}) = \frac{b_n^2}{2n}E(f^2(X)) + O\left(\frac{b_n^3}{n^2}\right)$$

as $n \to \infty$. Hence, (2.9) and (2.10) hold. $\square$

LEMMA 3. *Let $D$ satisfy (1.6) and assume $a_0$ is the dominating point for $(D, \gamma)$ with $g = t_0 f \in B^*$ as in (1.7) and (1.8). Then*

$$P(S_n/b_n \in D)$$
$$= \exp\{-b_n^2 n^{-1}\lambda_\gamma(a_0) - b_n^2 n^{-1}[\log \hat{\gamma}(g) - nb_n^{-2} \log E(e^{g(b_n S_n/n)})]\}J_n,$$

*where*

$$J_n = E(\exp\{-b_n^2 n^{-1}g(S_n/b_n - a_0)\} \quad (2.11)$$
$$\times e^{g(b_n S_n/n)}I(S_n/b_n \in D))/E(e^{g(b_n S_n/n)}).$$

*Furthermore, if $b_n = o(n^{2/3})$, then*

$$P(S_n/b_n \in D) \sim \exp\{-b_n^2 n^{-1}\lambda_\gamma(a_0)\}J_n \quad (2.12)$$

*and*

$$J_n \sim E(e^{-g(T_n - E(T_n))}I(T_n \in b_n^2 D/n)), \quad (2.13)$$

*where $T_n = \frac{b_n}{n}(Z_{n,1} + \cdots + Z_{n,n})$ and $Z_{n,1}, Z_{n,2}, \ldots$ are i.i.d. copies of $Z^{(n)}$ as defined in (1.13).*

PROOF. The proof of the representation formula for $P(S_n/b_n \in D)$ and (2.11) is simple algebra once one takes into account (1.7)(iii). Furthermore, if $b_n = o(n^{2/3})$, then (2.10) with $f = g$ implies (2.12) since $\log \hat{\gamma}(g) = \frac{1}{2}\int_B g^2(x) d\gamma(x) = \frac{1}{2}\int_B g^2(x) d\mu(x)$.

To verify (2.13) we observe that

$$J_n = \int_{B^n} \exp\{-g(x_1 + \cdots + x_n)b_n/n + b_n^2 n^{-1} g(a_0)\}$$
$$\times I(x_1 + \cdots + x_n \in b_n D) \, d\rho(x_1) \cdots d\rho(x_n),$$

where $\rho = \mathcal{L}(Z^{(n)})$. Thus,

$$J_n = E(\exp\{-g(T_n) + b_n^2 n^{-1} g(a_0)\} I(T_n \in b_n^2 n^{-1} D))$$
$$= \exp\{g(b_n^2 n^{-1} a_0 - E(T_n))\} E(\exp\{-g(T_n - E(T_n))\} I(T_n \in b_n^2 n^{-1} D)).$$

Now (2.1) implies $E(T_n) = (b_n^2/n)a_0 + O(b_n^3/n^2)$, and, hence, if $b_n = o(n^{2/3})$, (2.13) holds. $\square$

*Combining the lemmas.* Since (1.14) follows from Lemma 1 and (2.12) and (2.13) of Lemma 3 imply (1.12), Theorem 1 will follow once (1.15) is shown to hold. Using (1.12) we will have (1.15) provided

$$(2.14) \quad \limsup_{n \to \infty} n^{-1/2} b_n E(e^{-g(T_n - E(T_n))} I(T_n \in b_n^2 D/n)) \leq (2\pi \sigma_g^2)^{-1/2}.$$

Now

$$I_n \equiv E\left[e^{-g(T_n - E(T_n))} I\left(T_n \in \frac{b_n^2 D}{n}\right)\right]$$
$$(2.15) \quad = E\left[e^{-g(T_n - E(T_n))} I\left(T_n - E(T_n) \in \frac{b_n^2}{n}(D - n b_n^{-2} E(T_n))\right)\right]$$
$$= E\left[e^{-g(T_n - E(T_n))} I\left(T_n - E(T_n) \in \frac{b_n^2}{n}(D - a_0) + O\left(\frac{b_n^3}{n^2}\right)\right)\right],$$

where the last equality follows from (1.14) and that $T_n = \frac{b_n}{n}(Z_{n,1} + \cdots + Z_{n,n})$. If $\sigma_{T_n}^2 = (b_n^2/n)\sigma_{g,n}^2$ denotes the variance of $g(T_n)$, then

$$I_n \leq E(e^{-g(T_n - E(T_n))} I(g(T_n - E(T_n))/\sigma_{T_n} \geq O(b_n^3/n^2)/\sigma_{T_n}),$$

since $g(x) \geq 0$ for all $x \in D - a_0$. Here the term $O(b_n^3/n^2)$ may be positive or negative and $\sigma_{T_n} \to \infty$. Therefore,

$$I_n \leq \int_{]-\alpha_n, \infty[} e^{-\sigma_{T_n} u} \, dF_n(u),$$

where $0 \leq \alpha_n = O(b_n^2/n^{3/2})$ and $F_n$ denotes the distribution function of $g(T_n - E(T_n))/\sigma_{T_n}$. Thus, with $\Phi(u)$ the distribution function of a standard normal random variable, we have

$$I_n \leq \int_{]-\alpha_n, \infty[} \int_{[u, \infty[} \sigma_{T_n} e^{-\sigma_{T_n} x} \, dx \, dF_n(u)$$
$$= \int_{]-\alpha_n, \infty[} (F_n(x) - F_n(-\alpha_n)) \sigma_{T_n} e^{-\sigma_{T_n} x} \, dx$$

$$\leq \int_{]-\alpha_n,\infty[} (\Phi(x) - \Phi(-\alpha_n))\sigma_{T_n} e^{-\sigma_{T_n} x}\,dx$$
$$+ C_{\mathrm{BE}} e^{\sigma_{T_n}\alpha_n} E(|g(Z^{(n)})|^3)/(\sqrt{n}\sigma_{g,n}^3)$$
$$= \int_{]-\alpha_n,\infty[} e^{-\sigma_{T_n} u}\,d\Phi(u) + C_{\mathrm{BE}} e^{\sigma_{T_n}\alpha_n} E(|g(Z^{(n)})|^3)/(\sqrt{n}\sigma_{g,n}^3)$$

by the Berry–Esseen theorem. Taking into account that $\sigma_{T_n} \to \infty$ and $\alpha_n \to 0$ as $n \to \infty$, it follows after an elementary calculation that, for large enough $n$,

$$\int_{-\alpha_n}^{\infty} e^{-\sigma_{T_n} u}\,d\Phi(u) = e^{\sigma_{T_n}^2/2}(1 - \Phi(\sigma_{T_n} - \alpha_n))$$
$$\leq e^{\sigma_{T_n}^2/2} e^{-(\sigma_{T_n}-\alpha_n)^2/2}/\{(\sigma_{T_n}-\alpha_n)\sqrt{2\pi}\}$$
$$\leq \frac{1}{\sqrt{2\pi}} e^{\sigma_{T_n}\alpha_n}/(\sigma_{T_n} - \alpha_n).$$

Recalling (1.1) and (2.2), we see that $\sigma_{T_n}\alpha_n \to 0$. Moreover, it follows that $\sigma_{T_n} \sim \sigma_g b_n/\sqrt{n}$ as $n \to \infty$. Thus, we have for all $\delta > 0$ if $n \geq n(\delta)$,

$$I_n \leq (2\pi\sigma_g^2)^{-1/2} n^{1/2} b_n^{-1}(1+\delta)[1+o(1)].$$

The last inequality follows since $E(|g(Z^{(n)})|^3) \sim E(|g(X)|^3)$, $\sigma_{g,n} \to \sigma_g$ and $\frac{b_n}{n} \to 0$. Since $\delta > 0$ is arbitrary, we have (2.14) and Theorem 1 is proved.

**3. Proof of Theorem 2.** Applying (1.12), relation (1.18) and, consequently, Theorem 2 will follow if we show

$$(3.0) \qquad \liminf_{n\to\infty} n^{-1/2} b_n E(e^{-g(T_n - E(T_n))} I(T_n \in b_n^2 D/n)) > 0.$$

Since $0 \notin \overline{D}$ and $D$ contains slices near $a_0$ whose diameters dominate $\tau(s) = \beta(s|\log s|)^{1/2}$ for $0 < s \leq \delta$, we have

$$(3.1) \quad M_s \cap (D - a_0) \supset \{y + sx_0 : y \in M_0, \|y\| \leq \tau(s)\}, \qquad 0 < s \leq \delta,$$

where

$$M_s = \{x : g(x) = sg(x_0)\}, \qquad x_0 \in B, g(x_0) > 0, \delta > 0, \beta > 0.$$

Here $g = t_0 f \in B^*$ is related to the dominating point $a_0$ of $D$ with respect to $\gamma$ as in (1.7) and (1.8). Thus, by rescaling (3.1) with $r = st, 0 < s \leq \delta$, we have

$$
\begin{aligned}
M_r \cap t(D - a_0) &\\
&= M_{st} \cap t(D - a_0) \\
(3.2) \qquad &= t(M_s \cap D - a_0) \supset \{t(y + sx_0) : y \in M_0, \|y\| \leq \tau(s)\} \\
&= \{w + rx_0 : w/t \in M_0, \|w/t\| \leq \tau(r/t)\} \\
&= \{w + rx_0 : w \in M_0, \|w\| \leq \beta t^{1/2}(r|\log r/t|)^{1/2}\}.
\end{aligned}
$$

Hence,

$$t(D - a_0) \supset \{x = w + rx_0 : w \in M_0, 0 < r \leq t\delta, \tag{3.3}$$
$$\|w\| \leq \beta t^{1/2}(r|\log r/t|)^{1/2}\}.$$

Setting $\pi_g(x) = g(x)/g(x_0)$, we see $x - \pi_g(x)x_0 \in M_0$, and, thus, (3.3) implies

$$t(D - a_0) \supset \{x = x - \pi_g(x)x_0 + \pi_g(x)x_0 : 0 < \pi_g(x) \leq t\delta, \tag{3.4}$$
$$\|x - \pi_g(x)x_0\| \leq \beta t^{1/2}(|\pi_g(x)\log(\pi_g(x)/t)|)^{1/2}\}.$$

Recall that $\tilde{T}_n = T_n - E(T_n)$. Now (1.14) implies that $E(T_n) = b_n^2 a_0/n + \lambda_n$, where $\|\lambda_n\| = O(b_n^3/n^2)$, and, therefore, we have

$$(3.5) E\left(e^{-g(T_n - E(T_n))} I\left(T_n \in \frac{b_n^2}{n} D\right)\right) = E\left(e^{-g(\tilde{T}_n)} I\left(\tilde{T}_n \in \frac{b_n^2}{n}(D - a_0) + \lambda_n\right)\right).$$

Then, for $0 < A < B$,

$$\theta_n = (|\pi_g(\tilde{T}_n - \lambda_n)\log(\pi_g(\tilde{T}_n - \lambda_n)n/b_n^2)|)^{1/2},$$

$$\psi_n = \left(\frac{b_n^2}{n}\right)^{1/2} \left(\frac{A}{2g(x_0)}\right)^{1/2} \left|\log \frac{2Bn}{g(x_0)b_n^2}\right|^{1/2}$$

and $t = b_n^2/n$ in (3.4) implies for $n$ sufficiently large that

$$e^{2B} E(e^{-g(\tilde{T}_n)} I(\tilde{T}_n \in n^{-1}b_n^2(D - a_0) + \lambda_n))$$

$$\geq P\left(\frac{A}{2} < g(\tilde{T}_n) < 2B, \tilde{T}_n - \lambda_n \in n^{-1}b_n^2(D - a_0)\right)$$

$$\geq P\left(\frac{A}{2} < g(\tilde{T}_n) < 2B, 0 < \pi_g(\tilde{T}_n - \lambda_n) \leq \frac{b_n^2\delta}{n},\right.$$

(3.6)
$$\left.\|(\tilde{T}_n - \lambda_n) - \pi_g(\tilde{T}_n - \lambda_n)x_0\| \leq \beta\left(\frac{b_n^2}{n}\right)^{1/2} \theta_n\right)$$

$$\geq P(A < g(\tilde{T}_n) < B, \|\tilde{T}_n - \pi_g(\tilde{T}_n)x_0\| \leq \beta\psi_n - \|\lambda_n - \pi_g(\lambda_n)x_0\|)$$

$$\geq P\left(A < g(\tilde{T}_n) < B, \|\tilde{T}_n - \pi_g(\tilde{T}_n)x_0\| \leq \frac{\beta}{2}\psi_n\right).$$

The third inequality in (3.6) requires $n$ sufficiently large so that $A < g(\tilde{T}_n) < B$ implies $A/(2g(x_0)) < \pi_g(\tilde{T}_n - \lambda_n) < 2B/g(x_0) \leq \frac{b_n^2 \delta}{n}$ and this is immediate since $\lambda_n \to 0$, $\pi_g(x) = \frac{g(x)}{g(x_0)}$, and $b_n^2/n \to \infty$. The last inequality requires $n$ sufficiently large so that

$$\|\lambda_n - \pi_g(\lambda_n)x_0\| \leq \frac{\beta}{2}\left(\frac{b_n^2}{n}\right)^{1/2} \left(\frac{A}{2g(x_0)}\right)^{1/2} \left|\log\left(\frac{2Bn}{g(x_0)b_n^2}\right)\right|^{1/2}$$

and this is trivial since $\|\lambda_n\| \to 0$ and $b_n^2/n \to \infty$. Thus, for $n$ sufficiently large,

$$
\begin{aligned}
&e^{2B} E(e^{-g(\tilde{T}_n)} I(\tilde{T}_n \in n^{-1} b_n^2 (D - a_0) + \lambda_n)) \\
&\qquad \geq P(A < g(\tilde{T}_n) < B) \\
&\qquad - P\bigg(\|\tilde{T}_n - \pi_g(\tilde{T}_n) x_0\| \\
&\qquad\qquad > \frac{\beta}{2}\bigg(\frac{b_n^2}{n}\bigg)^{1/2} \bigg(\frac{A}{2g(x_0)}\bigg)^{1/2} \bigg|\log \frac{2Bn}{g(x_0) b_n^2}\bigg|^{1/2}\bigg).
\end{aligned}
$$
(3.7)

Defining again $\sigma_g^2 = E(g^2(X))$ and $\sigma_{g,n}^2 = E(g^2(Z^{(n)} - E(Z^{(n)})))$, it is evident that $\sigma_n^2 \equiv \sigma_{g(\tilde{T}_n)}^2 = \frac{b_n^2}{n} \sigma_{g,n}^2$ and the Berry–Esseen theorem implies that uniformly in $u \leq v$ and $n \geq 1$,

$$
\begin{aligned}
&|P(u/\sigma_n < g(\tilde{T}_n)/\sigma_n < v/\sigma_n) - P(u/\sigma_n < G < v/\sigma_n)| \\
&\qquad \leq C_{\mathrm{BE}} E(|g(Z^{(n)})|^3)/(\sigma_{g,n}^3 \sqrt{n}),
\end{aligned}
$$

where $G$ is standard normal. Now $E(|g(Z^{(n)})|^3) \sim E(|g(X)|^3)$ and $\sigma_{g,n} \sim \sigma_g > 0$ as $n \to \infty$. We thus have for large $n$,

$$
\begin{aligned}
&P(u/\sigma_n < g(\tilde{T}_n)/\sigma_n < v/\sigma_n) \\
&\qquad \geq P(u/\sigma_n < G < v/\sigma_n) - 2 C_{\mathrm{BE}} E(|g(X)|^3)/(\sigma_g^3 \sqrt{n}).
\end{aligned}
$$

Since $\sigma_n \to \infty$, we have

$$
P(u/\sigma_n < G < v/\sigma_n) \sim (v - u)/(2\pi \sigma_n^2)^{1/2}
$$

and, therefore, if $(v - u)/(2\pi)^{1/2} > 4 C_{\mathrm{BE}} E(|g(X)|^3)/\sigma_g^3$ we have

$$
P(u/\sigma_n < G < v/\sigma_n) \geq (1/2)(v - u)/(2\pi b_n^2 \sigma_g^2/n)^{1/2}
$$

because $b_n^2/n < n$. Taking $A = u$, $B = v$, we have

(3.8) $\qquad P(A < g(\tilde{T}_n) < B) \geq (1/2)(B - A) n^{1/2}/(2\pi \sigma_g^2 b_n^2)^{1/2},$

for all $n$ sufficiently large.

We now need an upper bound for

$$
P(\|\tilde{T}_n - \pi_g(\tilde{T}_n) x_0\| > (\beta/2) \psi_n) \leq P(\|\tilde{T}_n\| \geq (\beta/2K) \psi_n),
$$

where $K = \|Q\| < \infty$ and $Q: B \to B$ is the continuous operator given by $Q(x) = x - \pi_g(x) x_0$, $x \in B$.

To that end we first derive an upper bound for $E(\|\tilde{T}_n\|)$ where the following lemma comes in handy.

LEMMA 4. *Let $Y_1, \ldots, Y_n$ be i.i.d. random variables. Assume that*

$$P\left\{\left\|\sum_{j=1}^{n} Y_j\right\| \geq t_0\right\} \leq 10^{-4}.$$

*Then we have*

$$E\left(\left\|\sum_{j=1}^{n} Y_j\right\|\right) \leq 122 E\left(\max_{1 \leq j \leq n} \|Y_j\|\right) + 10^4 t_0.$$

PROOF. Using inequality (1.2.4) on page 10 in [10] with $s = t = u$, it follows that

$$P\left\{\left\|\sum_{j=1}^{n} Y_j\right\| > 61s\right\} \leq P\left\{\max_{1 \leq j \leq n} \|Y_j\| > s\right\} + 81\left(P\left\{\left\|\sum_{j=1}^{n} Y_j\right\| > s\right\}\right)^2$$

from which the moment inequality readily follows after integration by parts. □

Since $\tilde{T}_n = T_n - E(T_n)$ with $T_n = \frac{b_n}{n}(Z_{n,1} + \cdots + Z_{n,n})$, we have that $\{\tilde{T}_n/(b_n^2/n)^{1/2}\}$ is bounded in probability and Lemma 4, in conjunction with the Hölder inequality, implies for some $\alpha > 0$,

(3.9)  $$\limsup_{n \to \infty} E\|\tilde{T}_n\|/(b_n^2/n)^{1/2} \leq 10^4 \alpha,$$

using that $E(\max_{1 \leq j \leq n} \|Z_{n,j}\|^3)^{1/3} \leq n^{1/3} E(\|Z^{(n)}\|^3)^{1/3} \sim n^{1/3} E(\|X\|^3)^{1/3}$.

Thus, the Fuk–Nagaev inequality as given in [11], page 338, and that $|\log(n/b_n^2)| \to \infty$ implies

$$P(\|\tilde{T}_n\| > (\beta/2K)\psi_n)$$
$$\leq 9 \cdot 2^{11} t^{-3} E(\|Z^{(n)}\|^3)/n^{1/2} + \exp\{-t^2/(96 E\|Z^{(n)}\|^2)\},$$

where $t = (\beta/2K)(A/(2g(x_0)))^{1/2}|\log(2Bn/g(x_0)b_n^2)|^{1/2}$.

Since $E\|Z^{(n)}\|^3 \to E\|X\|^3$ and $E\|Z^{(n)}\|^2 \to E\|X\|^2$, we see that by taking $B = 2A$ and $A$ sufficiently large so that $\beta^2 A/(8g(x_0)) > 192 K^2 E\|X\|^2$, then this last probability is $o((b_n/n^{1/2})^{-1})$ as $n \to \infty$. Recalling (3.7) and (3.8), we can conclude that

$$e^{2B} E(e^{-g(\tilde{T}_n)} I(\tilde{T}_n \in n^{-1} b_n^2 (D - a_0) + \lambda_n)) \geq A/(4(2\pi \sigma_g^2 b_n^2/n)^{1/2})$$

for $n$ sufficiently large. Thus, (3.0) holds and Theorem 2 is established.

**4. Proof of Theorem 3.** Let $a_0$ be the unique dominating point of $(D, \gamma)$ and $g = t_0 f \in B^*$ be related to $a_0$ as in (1.7) and (1.8). Let $T_n = \frac{b_n}{n}(Z_{n,1} + \cdots + Z_{n,n})$ and $\tilde{T}_n = T_n - E(T_n)$ as before. As in the previous section we have to prove that

$$\liminf_{n \to \infty} n^{-1/2} b_n E(e^{-g(\tilde{T}_n)} I(T_n \in b_n^2 D/n)) > 0. \tag{4.1}$$

Under the present assumption on the set $D$ we obtain by the same argument as in Section 3 that

$$\begin{aligned}t(D - a_0) \supset \{x = x - \pi_g(x)x_0 + \pi_g(x)x_0 : \\ 0 < \pi_g(x) \le t\delta, \|x - \pi_g(x)x_0\| \le \beta t^{1/2} |\pi_g(x)|^{1/2}\}.\end{aligned} \tag{4.2}$$

Using again the fact that $ET_n = b_n^2 a_0/n + \lambda_n$, where $\lambda_n \in H, \lambda_n \to 0$, we have for any $A > 0$ that

$$\begin{aligned}E(e^{-g(\tilde{T}_n)} I(T_n \in b_n^2 D/n)) \\ = E(e^{-g(\tilde{T}_n)} I(\tilde{T}_n \in n^{-1} b_n^2 (D - a_0) + \lambda_n)) \\ \ge e^{-2A} P\{A < g(\tilde{T}_n) < 2A, 0 < \pi_g(\tilde{T}_n - \lambda_n) \le b_n^2 \delta/n, \\ \|(\tilde{T}_n - \lambda_n) - \pi_g(\tilde{T}_n - \lambda_n)x_0\| \\ \le \beta(b_n/n^{1/2})|\pi_g(\tilde{T}_n - \lambda_n)|^{1/2}\},\end{aligned} \tag{4.3}$$

which is for $n$ sufficiently large, greater than or equal to

$$\begin{aligned}e^{-2A} P\{A < g(\tilde{T}_n) < 2A, \\ \|(\tilde{T}_n - \lambda_n) - \pi_g(\tilde{T}_n - \lambda_n)x_0\| \le \beta(b_n/n^{1/2})|\pi_g(\tilde{T}_n - \lambda_n)|^{1/2}\}.\end{aligned}$$

This follows since $\pi_g(\lambda_n) \to 0$ and $b_n^2/n \to \infty$ imply eventually

$$\{A < g(\tilde{T}_n) < 2A\} \subset \{0 < \pi_g(\tilde{T}_n - \lambda_n) \le b_n^2 \delta/n\}.$$

Next, observe that also eventually

$$\{g(\tilde{T}_n) < 2A\} \subset \{\|x_0\| |\pi_g(\tilde{T}_n)|^{1/2} \le (\beta/4) b_n/n^{1/2}\},$$

which along with the fact that $\|\lambda_n\| \to 0$ implies for large $n$

$$\begin{aligned}P\{A < g(\tilde{T}_n) < 2A, \\ \|(\tilde{T}_n - \lambda_n) - \pi_g(\tilde{T}_n - \lambda_n)x_0\| \le \beta(b_n/n^{1/2})|\pi_g(\tilde{T}_n - \lambda_n)|^{1/2}\} \\ \ge P\{A < g(\tilde{T}_n) < 2A, \|\tilde{T}_n\| \le (\beta/2)(b_n/n^{1/2})|\pi_g(\tilde{T}_n - \lambda_n)|^{1/2}\}.\end{aligned}$$

Moreover, we have on the event $\{g(\tilde{T}_n) > A\}$ eventually, $|\pi_g(\tilde{T}_n - \lambda_n)| \ge |\pi_g(\tilde{T}_n)|/2$, hence, the last probability is, for large $n$,

$$\ge P\{A < g(\tilde{T}_n) < 2A, \|\tilde{T}_n\| \le (\beta/3)(b_n/n^{1/2})|\pi_g(\tilde{T}_n)|^{1/2}\}.$$

Recalling (4.3), we see that for large enough $n$,

(4.4)
$$\begin{aligned}E(e^{-g(\tilde{T}_n)}&I(T_n \in b_n^2 D/n))\\ &\geq e^{-2A}[P\{A < g(\tilde{T}_n) < 2A\}\\ &\quad - P\{\|\tilde{T}_n\| > (\beta/3)(b_n/n^{1/2})|\pi_g(\tilde{T}_n)|^{1/2}, \pi_g(\tilde{T}_n) \geq 0\}].\end{aligned}$$

In view of (3.8) we have if $A > 4C_{\mathrm{BE}}\sqrt{2\pi}E(|g(X)|^3)/\sigma_g^3$ for large $n$,

$$P(A < g(\tilde{T}_n) < 2A) \geq (1/2)An^{1/2}/(2\pi\sigma_g^2 b_n^2)^{1/2}.$$

Hence, by taking $A$ sufficiently large we will have (4.1), provided we show

(4.5) $$\limsup_{n\to\infty} \frac{b_n}{n^{1/2}} P(\|\tilde{T}_n\|^2 > \frac{\beta^2}{9}(n^{-1}b_n^2 \pi_g(\tilde{T}_n)), \pi_g(\tilde{T}_n) \geq 0) < \infty.$$

This will follow from the proposition below. Therefore, by combining (4.4) and (4.5) we have (4.1), and Theorem 3 is proved. □

**5. Proof of (4.5).** We will obtain a slightly more general result than needed.

PROPOSITION 1. *Let $X_{n,1},\ldots,X_{n,n}, n \geq 1$, be a triangular array of rowwise i.i.d. random vectors with values in the Hilbert space $H$ such that $E(X_{n,1}) = 0$ and $\sup_{n\geq 1} E(\|X_{n,1}\|^3) \leq M$. Let $\rho_n \to \infty$ such that $\rho_n = O(n^{1/2})$, $\|a\| = 1$, $f(x) = \langle x, a \rangle$, and assume that $\inf_{n\geq 1} E(f^2(X_{n,1})) \geq \delta^2 > 0$. If $S_n = \sum_{i=1}^n X_{n,i}$, then*

(5.1) $$\varlimsup_n \rho_n P(\|S_n/n^{1/2}\|^2 > \rho_n f(S_n/n^{1/2}), f(S_n) \geq 0) < \infty.$$

REMARK 1. It is easily checked that we can apply the above proposition with $X_{n,i} = Z_{n,i} - EZ_{n,i}$, $1 \leq i \leq n$, $n \geq 1$, so that this result indeed implies (4.4). [Recall that $\tilde{T}_n = (b_n/n)\sum_{i=1}^n (Z_{n,i} - EZ_{n,i})$.] The linear functional $\pi_g(x) = g(x)/g(x_0)$ can, of course, be normalized to have norm one without loss of generality.

REMARK 2. In the special case $\rho_n = n^{1/2}$, Proposition 1 above follows from Proposition 1 in [12] since

(5.2)
$$\begin{aligned}P(\|S_n/n^{1/2}\|^2 &> \rho_n f(S_n/n^{1/2}), f(S_n) \geq 0)\\ &= P(f(S_n) \geq 0) - P(\|S_n/n^{1/2}\|^2 \leq \rho_n f(S_n/n^{1/2}), f(S_n) \geq 0).\end{aligned}$$

REMARK 3. If $\rho_n = O(n^{1/2}/(\log n)^3)$, Proposition 1 also follows from Proposition 2 below (see Remark 6). Given that we consider in this paper only sequences $\rho_n$ of order $o(n^{1/6})$, this is more than sufficient for the proof of Theorem 3. We chose to include the proof via $U$-statistics as it allows a slightly larger $\rho_n$ which may be of future use.

PROOF OF PROPOSITION 1. In view of (5.2) it suffices to show that under the assumptions of Proposition 1 we have

$$|P(\|S_n/n^{1/2}\|^2 \leq \rho_n f(S_n/n^{1/2})) - 1/2| = O(\rho_n^{-1}).$$

This follows by applying the version of Theorem 1 of Alberink appearing on page 522 of [1]. Applying this result exactly as in the proof of Proposition 1 in [12], one obtains after some obvious modifications Proposition 1 above. □

**6. Proof of Theorem 4.** We prove this result for $d \geq 2$ only, though our proof can be modified to include the case $d = 1$ as well. However, in this case the result is well known and it can be proved more directly.

First observe that $H_\mu = \mathbb{R}^d$ and (1.19) implies that $E(e^{t\|X\|}) < \infty$ for some $t > 0$, where $\|\cdot\|$ is the usual Euclidean norm on $\mathbb{R}^d$. Hence, all possible moments of $X$ are finite, and Theorem 1 implies

$$(6.1) \qquad P(S_n/b_n \in D) \sim \exp\{-n^{-1}b_n^2 \lambda_\gamma(a_0)\} I_n,$$

where

$$(6.2) \qquad I_n = E(e^{-g(T_n - E(T_n))} I(T_n \in b_n^2 D/n)).$$

Recalling (1.14) and that $T_n = \frac{b_n}{n}(Z_{n,1} + \cdots + Z_{n,n})$, we also have

$$(6.3) \quad I_n = E(e^{-g(T_n - E(T_n))} I(T_n - E(T_n) \in b_n^2 n^{-1}(D - a_0) + \alpha_n)),$$

where $\alpha_n$ is a deterministic vector such that $\alpha_n = b_n^2 n^{-1} a_0 - E(T_n)$ and $\|\alpha_n\| = O(b_n^3/n^2) = o(1)$. Now $n^{1/2} b_n^{-1}(T_n - E(T_n)) = \sum_{j=1}^n (Z_{n,j} - E(Z_{n,j}))/n^{1/2}$ and, hence, if $G'_n$ is a mean zero Gaussian random vector with values in $\mathbb{R}^d$ and $\text{cov}(G'_n) = \text{cov}(Z^{(n)})$, then $n^{1/2} b_n^{-1}(T_n - E(T_n))$ can be approximated by $G'_n$. In particular, if we use the main result of Zaitsev [20] we have if $n$ is large enough for $\varepsilon > 0$ and all Borel-subsets $A$ of $\mathbb{R}^d$,

$$(6.4) \; P(n^{1/2} b_n^{-1}(T_n - E(T_n)) \in A) \leq P(G'_n \in A^\varepsilon) + c_1 \exp(-c_2 n^{1/2} \varepsilon/\tau),$$

where as usual $A^\varepsilon = \{x \in \mathbb{R}^d : \exists y \in A : \|x - y\| < \varepsilon\}$. Here $c_1, c_2$ are positive constants depending on $d$, and $\tau > 0$ depends on the distribution of $X$. To see this we note that from $E(e^{\|X\|/\tau}) < \infty$ for $\tau$ sufficiently large and $d\mathcal{L}(Z^{(n)})/d\mu(x) = e^{g(b_n x/n) - \log \hat{\mu}(b_n g/n)}$ with $b_n/n \to 0$, it follows that the distributions of $Z^{(n)}$ satisfy the hypothesis of Theorem 1.1 of Zaitsev [20] for $n \geq n_0$ and $\tau$ sufficiently large. This requires an elementary argument which we leave for the reader.

Hence, if we assume that the underlying $p$-space $(\Omega, \mathcal{F}, P)$ is rich enough, we can infer via the Strassen–Dudley theorem that for large enough $n$ and any given $\varepsilon > 0$, one can construct a mean zero Gaussian random vector $G'_{n,\varepsilon}$ with the same distribution as $G'_n$ so that

$$(6.5) \quad P(\|G'_{n,\varepsilon} - n^{1/2} b_n^{-1}(T_n - E(T_n))\| \geq \varepsilon) \leq c_1 \exp(-c_2 n^{1/2} \varepsilon/\tau).$$

To simplify our notation we set $\rho_n = b_n/n^{1/2}$. Choosing $\varepsilon = \varepsilon_n = \frac{1}{2}\rho_n^{-3/2}$ and writing $G'_n$ instead of $G'_{n,\varepsilon_n}$, we thus have if $n$ is large enough,

$$(6.6) \quad P(\|G'_n - \rho_n^{-1}(T_n - E(T_n))\| \geq \rho_n^{-3/2}/2)$$
$$\leq c_1 \exp(-c_2 n^{1/2} \rho_n^{-3/2}/2\tau) = o(n^{-1}).$$

We furthermore can assume that $G'_n = B_n Z$, where $Z$ is normal$(0, I)$-distributed and $B_n$ is a positive semi-definite, symmetric matrix so that $B_n^2 = \text{cov}(G'_n)$. ($I$ is the identity matrix.)

Set $G = BZ$, where $B$ is a positive definite, symmetric matrix so that $B^2 = \text{cov}(X)$ and $Z$ is as above. Arguing as in the proof of (2.2), we find that

$$(6.7) \quad \|B_n^2 - B^2\| = \|\text{cov}(Z^{(n)}) - \text{cov}(X)\| = O(b_n/n),$$

where $\|D\| = \sup_{\|x\| \leq 1} \|Dx\| = \sup_{\|x\| \leq 1} |\langle x, Dx \rangle|$ for symmetric $(d, d)$-matrices $D$.

Using the fact that $B$ is positive definite, one can infer (see Lemma 11) that

$$(6.8) \quad \|B_n - B\| = O(b_n/n),$$

which in turn via a standard exponential inequality for normal random vectors implies

$$(6.9) \quad P(\|G - G'_n\| \geq \rho_n^{-3/2}/2) \leq P(\|Z\| \geq \rho_n^{-3/2}/(2\|B_n - B\|)) = o(n^{-1}),$$

and we can conclude that

$$(6.10) \quad P(\|G - \rho_n^{-1}(T_n - E(T_n))\| \geq \rho_n^{-3/2}) = o(n^{-1}).$$

Set $C = D - a_0$. Returning to the integral (6.3), we can now infer that

$$(6.11) \quad I_n \leq I'_n + I''_n,$$

where

$$I'_n = E(e^{-g(T_n - E(T_n))} I(\rho_n^{-1}(T_n - E(T_n))) \in \rho_n C + \alpha_n/\rho_n,$$
$$\|\rho_n G - T_n + E(T_n)\| \leq \rho_n^{-1/2})$$

and

$$I''_n = E(e^{-g(T_n - E(T_n))} I(\rho_n^{-1}(T_n - E(T_n))) \in \rho_n C + \alpha_n/\rho_n,$$
$$\|\rho_n G - T_n + E(T_n)\| > \rho_n^{-1/2}).$$

Using the fact that $g(u) \geq 0, u \in C$, we readily obtain from (6.6) that

$$I''_n \leq e^{|g(\alpha_n)|} o(n^{-1}) = o(n^{-1}).$$

On the other hand, we have

$$I'_n \leq \exp(\|g\| \rho_n^{-1/2}) E(e^{-\rho_n g(G)} I(G \in (\rho_n C)^{\delta_n})),$$

where $\delta_n = \|\alpha_n\|\rho_n^{-1} + \rho_n^{-3/2} = o(\rho_n^{-1})$. As $g \geq -\|g\|\delta_n$ on $(\rho_n C)^{\delta_n}$, it easily follows from the subsequent Lemma 10 that

$$E(e^{-\rho_n g(G)} I(G \in (\rho_n C)^{\delta_n})) \leq E(e^{-\rho_n g(G)} I(G \in \rho_n C)) + e^{\|g\|\delta_n} O(\delta_n),$$

which in combination with the above estimates implies that

(6.12) $\quad I_n \leq \exp(\|g\|\rho_n^{-1/2}) E(e^{-\rho_n g(G)} I(G \in \rho_n C)) + o(\rho_n^{-1}).$

Changing in the proof of (6.12) the roles of $G$ and $\rho_n^{-1}(T_n - E(T_n))$ and setting $\alpha_n = 0$, we similarly get for large $n$,

(6.13) $\quad E(e^{-\rho_n g(G)} I(G \in \rho_n C)) \leq \exp(\|g\|\rho_n^{-1/2}) I_n + o(\rho_n^{-1}).$

More precisely, note that (6.10) implies that

$$P(\rho_n^{-1}(T_n - E(T_n)) \in (\rho_n C)^{\delta_n}) \leq P(G \in (\rho_n C)^{2\delta_n}) + o(n^{-1}),$$

which is on account of Lemma 10 and by a second application of (6.10), less than or equal to

$$P(G \in (\rho_n C)^{-2\delta_n}) + O(\delta_n) \leq P(\rho_n^{-1}(T_n - E(T_n)) \in \rho_n C) + O(\delta_n).$$

As Theorem 3 implies that $\liminf_{n \to \infty} I_n \rho_n > 0$, it is now evident that as $n \to \infty$,

(6.14) $\quad\quad\quad\quad\quad I_n \sim E(e^{-\rho_n g(G)} I(G \in \rho_n C)),$

where $G$ is a mean zero Gaussian random vector with covariance equal to that of $X$.

By the Cameron–Martin formula we have

(6.15) $\quad P\{G \in \rho_n D\} = \exp(-n^{-1} b_n^2 \lambda_\gamma(a_0)) E(e^{-\rho_n g(G)} I(G \in \rho_n C)),$

which in combination with (6.1) and (6.14) implies Theorem 4.

LEMMA 10. *Let $G$ be a centered, $\mathbb{R}^d$-valued, Gaussian random vector with covariance $V$ and support all of $\mathbb{R}^d$. If $\lambda$ is the minimal eigenvalue of $V$, then for all $\varepsilon > 0$ and all convex sets $C$, there exists a constant $c_d$ depending only on $d$ such that*

(6.16) $\quad\quad\quad\quad\quad P(G \in C^\varepsilon \setminus C^{-\varepsilon}) \leq 2c_d \lambda^{-1/2} \varepsilon,$

*where $C^\varepsilon = \bigcup_{x \in C} B(x, \varepsilon)$ and $C^{-\varepsilon} = \{x : B(x, \varepsilon) \subseteq C\}$.*

PROOF. If the covariance matrix $V$ is the identity matrix $I$, then this follows from Corollary 3.2 in [5] with $\lambda = 1$. Otherwise, let $A$ be a symmetric positive definite matrix such that $A^2 = V^{-1}$. Then $Z = AG$ has covariance $I$ on $\mathbb{R}^d$, and since $A$ has full rank,

$$P(G \in C^\varepsilon \setminus C^{-\varepsilon}) = P(Z \in T_A(C^\varepsilon) \setminus T_A(C^{-\varepsilon})),$$

where $T_A : \mathbb{R}^d \to \mathbb{R}^d$ is the linear operator determined by $A$. Noting that by an elementary argument,

$$T_A(C^\varepsilon) \subseteq T_A(C)^{\lambda^{-1/2}\varepsilon} \tag{6.17}$$

and

$$T_A(C^{-\varepsilon}) \supseteq (T_A(C))^{-\lambda^{-1/2}\varepsilon}, \tag{6.18}$$

we have by Corollary 3.2 in [5] that

$$P(G \in C^\varepsilon \setminus C^{-\varepsilon}) \leq P(Z \in (T_A(C))^{\lambda^{-1/2}\varepsilon} \setminus (T_A(C))^{-\lambda^{-1/2}\varepsilon})$$
$$\leq 2c_d \lambda^{-1/2}\varepsilon.$$

Hence, the lemma follows. □

We finally state a lemma from linear algebra which was needed for the above proof.

LEMMA 11. *Let $A, E$ be symmetric $(d,d)$-matrixes so that $A$ and $A + E$ are positive definite. Then we have for the positive definite square root matrices $\sqrt{A}$ and $\sqrt{A+E}$,*

$$\|\sqrt{A+E} - \sqrt{A}\| \leq C\|E\|,$$

*where $C$ is a positive constant depending on the smallest eigenvalue of $A$ and $A + E$.*

PROOF. The lemma is very easy to prove if $AE = EA$. In general, it follows from relation (X.46) on page 305 of [4] setting $r = 1/2$. □

**7. Proof of Theorem 5.** We still need the following lemma.

LEMMA 12. *Let $G$ be a centered Gaussian random variable on a separable Hilbert space $H$ and $D = \{x : \|x - a\| < R\}$, where $0 < R < \|a\|$, is an open ball in $H$ satisfying (1.6.ii) and (1.6.iii). Assume that $a_0 \in \partial D$ is the unique dominating point for $D$ with respect to $\gamma (=$ distribution of $G)$ and let $g$ be as in (1.7) and (1.8). Then we have the following for any positive sequence $\{b_n\}$ satisfying (1.1) and $\rho_n = b_n/n^{1/2}$:*

(i) $P(G \in \rho_n D) = \exp(-\rho_n^2 \lambda_\gamma(a_0)) E(e^{-\rho_n g(G)} I(G \in \rho_n(D - a_0)))$ *and*

(ii) $E(e^{-\rho_n g(G)} I(G \in \rho_n(D - a_0))) \sim \int_0^\infty e^{-s} P(\|G_2\|^2 \leq 2sbR^2) \, ds / (2\pi \sigma_g^2 \rho_n^2)^{1/2}$,

*as $n \to \infty$, where $\sigma_g^2 = E(g^2(G)), G_2 = G - G_1, b = 1/g(a - a_0)$ and $G_1 = g(G)E(Gg(G))/\sigma_g^2$.*

(iii) If $G_n$ is centered Gaussian with $\text{cov}(G_n) = \text{cov}(Z^{(n)})$, where $Z^{(n)}$ is defined as in Theorem 1, then

$$E(e^{-\rho_n g(G_n)} I(G_n \in \rho_n(D - a_0))) \tag{7.1}$$
$$\sim \int_0^\infty e^{-s} P(\|G_2\|^2 \leq 2sbR^2)\, ds / (2\pi \sigma_g^2 \rho_n^2)^{1/2}$$

as $n \to \infty$, where $G_1, G_2, \sigma_g^2$ and $b$ are as in (ii).

PROOF.  Part (i) follows directly from the definition "dominating point" and a simplification of the representation formula when $\mu$ is centered Gaussian. A key fact is that in this special case the law of $Z^{(n)}$ is that of $G + b_n a_0 / n$. This follows from the Cameron–Martin formula by an argument as in (2.5).

The proof of (ii) will follow along lines similar to those for (iii), so we now turn to the proof of (iii).

The proof of (iii) is as follows. Recall $\sigma_{g,n}^2 = E(g^2(G_n)) = E(g^2(Z^{(n)})) = \sigma_g^2 + O(b_n/n)$, and write $G_n = G_{n,1} + G_{n,2}$, where $G_{n,1} = g(G_n) E(G_n g(G_n))/\sigma_{g,n}^2$ and $G_{n,2} = G_n - G_{n,1}$.

Note that $g(G_{n,1}) = g(G_n)$, so $G_{n,2}$ has support in $\{x : g(x) = 0\}$. Furthermore, $G_{n,1}$ and $G_{n,2}$ are independent Gaussian random vectors and, hence, if

$$I_n(G) = E(e^{-\rho_n g(G)} I(G \in \rho_n(D - a_0))),$$

then for all $n$ sufficiently large,

$$I_n(G_n) = (2\pi \sigma_{g,n}^2)^{-1/2} \int_0^\infty e^{-\rho_n u} \tilde{h}(n, u) \exp\{-u^2 / 2\sigma_{g,n}^2\}\, du,$$

where

$$\tilde{h}(n, u) = P(G_{n,2} \in \rho_n(D - a_0) - u E(G_{n,1} g(G_{n,1}))/\sigma_{g,n}^2 | g(G_{n,1}) = u).$$

Thus, for sufficiently large $n$,

$$I_n(G_n) = (2\pi \rho_n^2 \sigma_{g,n}^2)^{-1/2} \int_0^\infty e^{-s} h(n, s) \exp\{-s^2/(2\rho_n^2 \sigma_{g,n}^2)\}\, ds \tag{7.2}$$

where $s = \rho_n u$, and since $G_{n,2}$ and $G_{n,1}$ are independent,

$$h(n, s) = P(G_{n,2} \in \rho_n(D - a_0) - \rho_n^{-1} s E(G_{n,1} g(G_{n,1}))/\sigma_{g,n}^2). \tag{7.3}$$

Now $D - a_0 = \{x : \|x - x_0\| < R\}$, where $x_0 = a - a_0$, and if $g(x_0) = 1/b$, we see that $bx_0 - E(G_{n,1} g(G_{n,1}))/\sigma_{g,n}^2$ is in $\{x : g(x) = 0\}$.

Furthermore, $\{x : g(x) = 0\}$ is tangent to the sphere $D - a_0$ at the origin and, hence, $x_0$ is perpendicular to the hyperplane $\{x : g(x) = 0\}$ as $D$ is a ball in Hilbert space. Thus, by the Pythogorean theorem, if $g(x) = 0$, then

$$x \in k(D - a_0) - bsx_0/k \quad \text{iff } \|x\|^2 < (kR)^2 - \left(k - \left(\frac{s}{k}\right)b\right)^2 R^2,$$

$$\text{iff } \|x\|^2 < 2sbR^2 - R^2b^2s^2/k^2.$$

Setting $E_n = \rho_n(D - a_0)$ and $k = \rho_n$ in the above, we therefore have

$$h(n, s) = P(G_{n,2} \in E_n - \rho_n^{-1} sbx_0 + \rho_n^{-1} s(bx_0 - E(G_{n,1}g(G_{n,1}))/\sigma_{g,n}^2))$$
$$= P(\|G_{n,2} - \rho_n^{-1} s(bx_0 - E(G_{n,1}g(G_{n,1}))/\sigma_{g,n}^2)\|^2$$
$$\leq 2sbR^2 - R^2b^2s\rho_n^{-2}).$$

Using the continuity of the distribution of the norm of a Gaussian random vector in a separable Hilbert space, and that $G_{n,2}$ converges weakly to $G_2$ on $\{x : g(x) = 0\}$, we thus see that

$$(7.4) \qquad \lim_{n \to \infty} h(n, s) = P(\|G_2\|^2 \leq 2sbR^2)$$

for $0 < s < \infty$. Combining (7.2)–(7.4), we thus have (7.1) since $\lim_n \sigma_{g,n}^2 = \sigma_g^2$. Hence, part (iii) of Lemma 12 is proved.

To verify the same asymptotics for $I_n(G)$ is quite similar with $G_2$ and $G_1$ replacing $G_{n,2}$ and $G_{n,1}$ throughout the argument. Hence, Lemma 12 is proved. □

*Conclusion of the proof of Theorem* 5. We first prove that

$$(7.5) \qquad \limsup_{n \to \infty} P(S_n \in b_n D)/P(G \in \rho_n D) \leq 1.$$

If $H$ is finite-dimensional, this follows from Theorem 4 and the usual isometry between $H$ and $\mathbb{R}^d$.

If $H$ is infinite-dimensional, take $\{e_1, e_2, \ldots\}$ to be a complete orthonormal sequence for $H$ with $e_1 = v/\|v\|$, where $v$ is the unique vector in $H$ so that $g(\cdot) = \langle v, \cdot \rangle$, and $g$ is as in (1.7). Define the orthogonal projection

$$(7.6) \qquad \pi_d(x) = \sum_{j=1}^{d} \langle e_j, x \rangle e_j, d \geq 1.$$

Then $g(x) = g(\pi_d(x))$ for all $x \in H$ and $d \geq 1$.

Applying (1.12), we have

$$(7.7) \qquad P(S_n \in b_n D) \sim \exp(-\rho_n^2 \lambda_\gamma(a_0)) I_n,$$

where

$$(7.8) \qquad I_n = E(\exp(-g(T_n - E(T_n))I(T_n \in \rho_n^2 D))).$$

Since $\pi_d : H \to H$ satisfies $g(x) = g(\pi_d(x))$, we easily have

$$
\begin{aligned}
(7.9) \quad I_n &\leq E(\exp(-g(T_n - E(T_n))I(\pi_d(T_n) \in \rho_n^2 \pi_d(D)))) \\
&= E(\exp(-g(\pi_d(T_n)) - E(\pi_d(T_n)))I(\pi_d(T_n) \in \rho_n^2 \pi_d(D))) =: I_{n,d}.
\end{aligned}
$$

Now by the proof of Theorem 4 [which also applies to the finite-dimensional space $\pi_d(H)$ by isometry],

$$(7.10) \qquad I_{n,d} \sim I_n(\pi_d(G)),$$

where

$$(7.11) \quad I_n(\pi_d(G)) = E(\exp(-\rho_n \pi_d(G))I(\pi_d(G) \in \rho_n \pi_d(D - a_0)))$$

and $G$ is a mean zero Gaussian random vector with covariance equal to that of $X$. Crucial to this last claim is the fact that $\pi_d(x_0) = x_0$, $\pi_d$ is an orthogonal projection,

$$
\begin{aligned}
(7.12) \quad \pi_d(D - a_0) &= \{\pi_d(y) : \|\pi_d(y - x_0)\| < R\} \\
&= \{\pi_d(y) : \|\pi_d(y) - x_0\| < R\},
\end{aligned}
$$

where $x_0 = a - a_0$, and the Radon–Nikodym derivative of the law of $\pi_d(Z^{(n)})$ with respect to the law of $\pi_d(X)$ is the same as that of the law of $Z^{(n)}$ with respect to the law of $X$. [Note that $\pi_d(x_0) = x_0$, since $D$ being a ball in Hilbert space implies that $v = \lambda x_0$ for some $\lambda > 0$ as the hyperplane $\{x : g(x) = 0\}$ is tangent to $D - a_0$ at zero.]

Hence, by Lemma 12(ii) applied to $\pi_d(G)$ in the subspace $\pi_d(H)$, we have

$$(7.13) \quad I_n(\pi_d(G)) \sim (2\pi \sigma_g^2 \rho_n^2)^{-1/2} \int_0^\infty e^{-s} P(\|\pi_d(G_2)\|^2 \leq 2sbR^2) \, ds,$$

where $g(\pi_d(x_0)) = g(x_0) = 1/b$, $\sigma_g^2 = E(g^2(X)) = E(g^2(\pi_d(X)))$, and $G_2 = G - G_1$ is a centered Gaussian random vector on $H$ with

$$(7.14) \qquad G_1 = g(G)E(Gg(G))/\sigma_g^2.$$

[Note that $\pi_d(G_1) = g(G)E(\pi_d(G)g(G))/\sigma_g^2$.]

Thus, when $D$ is a ball as indicated, for all $d \geq 2$ we have by (7.7), (7.9), (7.10), (7.13) and Lemma 12 that

$$(7.15) \quad \limsup_{n \to \infty} \frac{P(S_n \in b_n D)}{P(G \in \rho_n D)} \leq \frac{\int_0^\infty e^{-s} P(\|\pi_d(G_2)\|^2 \leq 2sbR^2) \, ds}{\int_0^\infty e^{-s} P(\|G_2\|^2 \leq 2sbR^2) \, ds}.$$

Letting $d \to \infty$, it easily follows by the dominated convergence theorem that the right-hand side approaches 1, which implies (7.5).

It remains to be shown that

$$(7.16) \qquad \liminf_{n \to \infty} P(S_n \in b_n D)/P(G \in \rho_n D) \geq 1.$$

But this follows from (1.12) in combination with Lemma 12 and the following proposition applied when the law of $X_{n,1}$ is equal to the law of $Z^{(n)} - E(Z^{(n)})$. To be more specific, let $A_n = S_n/n^{1/2}$, with $S_n$ as in Proposition 2, $f(x) = \langle x, x_0 \rangle$, $x \in H$ and notice that then

$$\begin{aligned}
I_n &= E(e^{-\rho_n g(A_n)} I(A_n \in \rho_n(D - a_0) + \alpha_n/\rho_n)) \\
&= E(e^{-\rho_n g(A_n)} I(\|A_n - \alpha_n/\rho_n - \rho_n x_0\|^2 \leq \rho_n^2 \|x_0\|^2)) \\
&= E(e^{-\lambda \rho_n f(A_n)} I(\|A_n - \alpha_n/\rho_n\|^2 \leq 2\rho_n f(A_n - \alpha_n/\rho_n)))
\end{aligned}$$
(7.17)

because $g(x) = \langle x, v \rangle = \lambda \langle x, x_0 \rangle = \lambda f(x)$. Also, recall that $\alpha_n = (b_n^2/n)a_0 - E(T_n)$ satisfies $\|\alpha_n\| = o(1)$.

Similarly, it follows that

$$\begin{aligned}
I_n(G_n) &= E(e^{-\rho_n g(G_n)} I(G_n \in \rho_n(D - a_0))) \\
&= E(e^{-\lambda \rho_n f(G_n)} I(\|G_n\|^2 \leq 2\rho_n f(G_n))).
\end{aligned}$$
(7.18)

Hence, if $\text{cov}(G_n) = \text{cov}(X_{n,1})$, then by Lemma 12(ii) and (iii) we have $I_n(G_n) \sim I_n(G)$ and by Proposition 2 (applied with $a = x_0/\|x_0\|$ and $\rho_n$ replaced by $\|x_0\|\rho_n$) that $\liminf_n I_n/I_n(G_n) \geq 1$, so the end result is that $\liminf_n I_n/I_n(G) \geq 1$, which proves (7.16). Thus, Theorem 5 follows once Proposition 2 is proved.

PROPOSITION 2. *Let $X_{n,1}, \ldots, X_{n,n}, n \geq 1$ be a triangular array of row-wise i.i.d. random vectors with values in the Hilbert space $H$ such that $E(X_{n,1}) = 0$ and $\sup_{n \geq 1} E(\|X_{n,1}\|^3) \leq M$. Let $\lambda > 0$ be a constant and $\rho_n \to \infty$ such that $\rho_n = O(n^{1/2}/(\log n)^3)$, $f(x) = \langle x, a \rangle$, where $\|a\| = 1$ and assume that $\inf_{n \geq 1} E(f^2(X_{n,1})) \geq \delta^2 > 0$. If $z_n$ is a sequence in $H$ with $\|z_n\| = o(\rho_n^{-1})$ and $S_n = \sum_{i=1}^n X_{n,i}$, then*

$$\liminf_{n \to \infty} E[\exp(-\lambda \rho_n f(S_n/n^{1/2})) I\{\|S_n/n^{1/2} + z_n\|^2 \leq 2\rho_n f(S_n/n^{1/2} + z_n)\}]/J_n$$
$$\geq 1,$$

*where $J_n = E[\exp(-\lambda \rho_n f(Y_n)) I\{\|Y_n\|^2 \leq 2\rho_n f(Y_n)\}]$ and $Y_n$ is a Gaussian mean zero random vector with covariance equal to that of $X_{n,1}$.*

REMARK 4. It is also possible to prove that

$$\limsup_{n \to \infty} E[\exp(-\lambda \rho_n f(S_n/n^{1/2})) I\{\|S_n/n^{1/2} + z_n\|^2 \leq 2\rho_n f(S_n/n^{1/2} + z_n)\}]/J_n$$
$$\leq 1,$$

so that we actually have an asymptotic equivalence. We did not work out the details since for the upper bound part of Theorem 5, it seems much more efficient to use the projection method as in the first part of Section 7.

REMARK 5. Given a fixed sequence $\rho_n$, one can replace the third moment assumption by some uniformity condition on the moments of order $2 + \eta$, where $0 < \eta \leq 1$ has to be determined depending on $\rho_n$.

REMARK 6. The subsequent proof also works for $\lambda = 0$. Following the proof until the inequality after (8.49), one sees that

$$P\{\|S_n/n^{1/2} + z_n\|^2 \leq 2\rho_n f(S_n/n^{1/2} + z_n)\}$$
$$\geq P\{\|Q_n(Y_n')\|^2 \leq 2\rho_n f(Y_n')\} - o(\rho_n^{-1}),$$

where $\|Q_n(Y_n')\|^2$ and $f(Y_n')$ are independent and $Y_n'$ is a Gaussian mean zero random vector. Choosing $z_n = 0$ and replacing $\rho_n$ by $\rho_n/2$, one readily obtains via Lemma 13 and the Berry–Esseen inequality that $\limsup_{n\to\infty} \rho_n P\{\|S_n/n^{1/2}\|^2 > \rho_n f(S_n/n^{1/2}), f(S_n) \geq 0\} < \infty$ provided that $\rho_n = O(n^{1/2}/(\log n)^3)$.

The proof of Proposition 2 is quite long. So it might be useful to give first an outline of the basic steps. To simplify our notation let

$$I_n := E[\exp(-\lambda \rho_n f(S_n/n^{1/2})) I\{\|S_n/n^{1/2} + z_n\|^2 \leq 2\rho_n f(S_n/n^{1/2} + z_n)\}].$$

From the proof of Theorem 3 it follows that $I_n$ is of order $O(\rho_n^{-1})$ so that it is sufficient to derive lower bounds up to terms of order $o(\rho_n^{-1})$.

We first show in step (i) that

$$I_n \geq I_{n,1} = E[\exp(-\lambda \rho_n f(\tilde{S}_n/n^{1/2})) I\{\|\tilde{S}_n/n^{1/2}\|^2 \leq (2 - \varepsilon_{n,1}) \rho_n f(\tilde{S}_n/n^{1/2})\}]$$
$$+ o(\rho_n^{-1}),$$

where $\tilde{S}_n$ are sums of truncated, centered random variables $\tilde{X}_{n,i}$ and $\varepsilon_{n,1} \to 0$. Note that this also shows that we can discard the vectors $z_n$.

Then we choose in step (ii) vectors $w_n$ so that the variables $f(\tilde{X}_{n,i})$ and $Q_n(\tilde{X}_{n,i})$ are uncorrelated and we show that

$$I_{n,1} \geq I_{n,2} = E[\exp(-\lambda \rho_n f(\tilde{S}_n/n^{1/2}))$$
$$\times I\{\|Q_n(\tilde{S}_n/n^{1/2})\|^2 \leq (2 - \varepsilon_{n,2}) \rho_n f(\tilde{S}_n/n^{1/2})\}] + o(\rho_n^{-1}),$$

where $Q_n(x) = x - f(x) w_n, x \in H$ and $\varepsilon_{n,2} \to 0$.

In step (iii) we smooth the variables $f(\tilde{S}_n/n^{1/2})$ and $\|Q_n(\tilde{S}_n/n^{1/2})\|^2$ by adding small independent normal variables and we show that

$$I_{n,2} \geq I_{n,3} = E[\exp(-\lambda \rho_n W_n) I\{V_n \leq (2 - \varepsilon_{n,3}) \rho_n W_n, W_n > 0\}] + o(\rho_n^{-1}),$$

where $W_n$ and $V_n$ are the smoothed variables and $\varepsilon_{n,3} \to 0$.

In step (iv) we make the crucial transition to the Gaussian case. We show that we can replace the variables $W_n, V_n$ by smoothed versions $\bar{W}_n$ and $\bar{V}_n$ of $f(Y_n')$ and $\|Q_n(Y_n')\|^2$, respectively. That is, we prove that

$$I_{n,3} \geq I_{n,4} = E[\exp(-\lambda \rho_n \bar{W}_n) I\{\bar{V}_n \leq (2 - \varepsilon_{n,3}) \rho_n \bar{W}_n, \bar{W}_n > 0\}] + o(\rho_n^{-1}),$$

where $Y'_n$ is mean zero Gaussian with $\text{cov}(Y'_n) = \text{cov}(\tilde{X}_{n,1})$. The crucial result for proving this last inequality is a certain local limit theorem, Lemma 18.

The proof of this lemma can be found in part (v) of the proof. As already mentioned in the Introduction we use an adaptation of the characteristic function method for proving Berry–Esseen type results in Hilbert space. In particular, we use a modification of a symmetrization lemma of Götze [13] [see (8.39)].

In step (vi) we then show that we can remove the smoothing variables, that is, we prove that

$$I_{n,4} \geq I_{n,5} = E[\exp(-\lambda\rho_n f(Y'_n))I\{\|Q_n(Y'_n)\|^2 \leq (2 - \varepsilon_{n,3})\rho_n f(Y'_n)\}] + o(\rho_n^{-1}).$$

Here it is very helpful that the variables $f(Y'_n)$ and $\|Q_n(Y'_n)\|^2$ are independent due to the choice of $w_n$ in step (ii).

In the following step (vii) we remove the sequence $\varepsilon_{n,3}$, that is, we prove that

$$I_{n,5} \geq I_{n,6} = E[\exp(-\lambda\rho_n f(Y'_n))I\{\|Q_n(Y'_n)\|^2 \leq 2\rho_n f(Y'_n)\}] + o(\rho_n^{-1}).$$

In the final step (viii) we use independence and the inequality of Anderson to prove that

$$I_{n,6} \geq J_n + o(\rho_n^{-1}).$$

**8. Proof of Proposition 2.** (i) Let $X'_{n,i} = X_{n,i}I(\|X_{n,i}\| \leq \delta n^{1/2})$, $\tilde{X}_{n,i} = X'_{n,i} - EX'_{n,i}$, $1 \leq i \leq n$, $n \geq 1$, and denote the corresponding sums by $S'_n$ and $\tilde{S}_n$, respectively.

Then it is easy to see that

(8.1) $$I_n \geq I'_n - \exp(\lambda\rho_n\|z_n\|)P(S_n \neq S'_n),$$

where $I'_n = E[\exp(-\lambda\rho_n f(S'_n/n^{1/2}))I\{\|S'_n/n^{1/2} + z_n\|^2 \leq 2\rho_n f(S'_n/n^{1/2} + z_n)\}]$. [Note that $f(S'_n/n^{1/2} + z_n) \geq 0$ implies $f(S'_n/n^{1/2}) \geq -\|z_n\|$.]

We have trivially, by Markov's inequality,

$$P(S_n \neq S'_n) \leq nP(X_{n,1} \neq X'_{n,1}) \leq M\delta^{-3}/\sqrt{n}.$$

Next, set $z'_n = ES'_n/n^{1/2}$, $z''_n = z_n + z'_n$ and observe that

$$\|z'_n\| \leq n^{1/2}E(\|X_{n,1}\|I(\|X_{n,1}\| > \delta n^{1/2})) \leq M\delta^{-2}/\sqrt{n}.$$

We can then further conclude from $|f(x)| \leq \|x\|$ that

$$I'_n = E[\exp(-\lambda\rho_n f(\tilde{S}_n/n^{1/2} + z'_n))I\{\|\tilde{S}_n/n^{1/2} + z''_n\|^2 \leq 2\rho_n f(\tilde{S}_n/n^{1/2} + z''_n)\}]$$
$$\geq \exp(-\lambda\rho_n\|z'_n\|)I''_n,$$

where $I''_n = E[\exp(-\lambda\rho_n f(\tilde{S}_n/n^{1/2}))I\{\|\tilde{S}_n/n^{1/2} + z''_n\|^2 \leq 2\rho_n f(\tilde{S}_n/n^{1/2} + z''_n)\}]$.

Let $A_n$ be the event $\{\|\tilde{S}_n/n^{1/2} + z_n''\|^2 \leq 2\rho_n f(\tilde{S}_n/n^{1/2} + z_n'')\}$. Then we clearly have

$$A_n \supset \{\|\tilde{S}_n/n^{1/2}\|^2 + 2\|z_n''\|\|\tilde{S}_n/n^{1/2}\| + \|z_n''\|^2$$
$$\leq 2\rho_n(f(\tilde{S}_n/n^{1/2}) + f(z_n''))\} =: A_n'.$$

Let $B_n = \{f(\tilde{S}_n/n^{1/2}) \leq \|z_n''\|\varepsilon_n^{-1}\}$, where $\varepsilon_n \searrow 0$ will be specified later. Consider further the event $C_n = \{\|\tilde{S}_n/n^{1/2}\|^2 \leq (2-\varepsilon_n)(1+\varepsilon_n)^{-2}\rho_n f(\tilde{S}_n/n^{1/2})\}$. Note that we have on the event $C_n \cap B_n^c$,

$$\|\tilde{S}_n/n^{1/2}\| \geq f(\tilde{S}_n/n^{1/2}) > \|z_n''\|\varepsilon_n^{-1}$$

and, consequently,

$$\|\tilde{S}_n/n^{1/2}\|^2 + 2\|z_n''\|\|\tilde{S}_n/n^{1/2}\| + \|z_n''\|^2 \leq \|\tilde{S}_n/n^{1/2}\|^2(1+\varepsilon_n)^2.$$

Furthermore, we have on this event $|f(z_n'')| \leq \varepsilon_n f(\tilde{S}_n/n^{1/2})$ and, thus,

$$2\rho_n(f(\tilde{S}_n/n^{1/2}) + f(z_n'')) \geq (2-\varepsilon_n)\rho_n f(\tilde{S}_n/n^{1/2}).$$

We see that $C_n \cap B_n^c \subset A_n'$, which in turn implies $I_{A_n} \geq I_{C_n} - I_{C_n \cap B_n}$.

Using the elementary inequality $(2-\varepsilon_n)(1+\varepsilon_n)^{-2} \geq 2 - 5\varepsilon_n$, we find that

$$I_n'' \geq E[\exp(-\lambda\rho_n f(\tilde{S}_n/n^{1/2}))I\{\|\tilde{S}_n/n^{1/2}\|^2 \leq (2-5\varepsilon_n)\rho_n f(\tilde{S}_n/n^{1/2})\}]$$
$$- E[\exp(-\lambda\rho_n f(\tilde{S}_n/n^{1/2}))I_{C_n \cap B_n}].$$

Recalling that $f(\tilde{S}_n/n^{1/2}) \geq 0$ on $C_n$, we further have

$$E[\exp(-\lambda\rho_n f(\tilde{S}_n/n^{1/2}))I_{C_n \cap B_n}]$$
$$\leq P(C_n \cap B_n) \leq P\{0 \leq f(\tilde{S}_n/n^{1/2}) \leq \|z_n''\|\varepsilon_n^{-1}\} =: p_n.$$

We need an upper bound for $p_n$. To that end we first note that

$$\tilde{\sigma}_{f,n}^2 = E(f^2(\tilde{X}_{n,1})) \geq E(f^2(X_{n,1})) - 2E(f^2(X_{n,1})I(\|X_{n,1}\| \geq \delta n^{1/2}))$$
$$\geq \delta^2 - 2M\delta^{-1}/n^{1/2} \geq \delta^2/2,$$

provided that $4M\delta^{-3}/n^{1/2} \leq 1$. Using the Berry–Esseen inequality, it now follows that

$$p_n \leq \|z_n''\|\varepsilon_n^{-1}\delta^{-1} + C_{\text{BE}}E|f(\tilde{X}_{n,1})|^3(\delta^2/2)^{-3/2}/n^{1/2}.$$

Employing the inequalities

(8.2) $$E|f(\tilde{X}_{n,1})|^3 \leq E\|\tilde{X}_{n,1}\|^3 \leq 8M,$$

we find that

$$p_n \leq \|z_n''\|\varepsilon_n^{-1}\delta^{-1} + 8^{3/2}C_{\text{BE}}M\delta^{-3}/n^{1/2},$$

which is trivially true if $4M\delta^{-3}/n^{1/2} > 1$.

Thus, $p_n$ has the order $O(n^{-1/2} \vee \|z_n''\|\varepsilon_n^{-1})$, which is of order $o(\rho_n^{-1})$ if $\varepsilon_n$ converges slowly enough to 0. (For instance, we can set $\varepsilon_n = \|z_n''\|^{1/2}\rho_n^{1/2}$.)

(ii) Let $w_n = \tilde{\sigma}_{f,n}^{-2} E(\tilde{X}_{n,1} f(\tilde{X}_{n,1}))$ and $Q_n(x) = x - f(x)w_n$ for $x \in H$. If $Y_n'$ is a mean zero Gaussian random vector with $\text{cov}(Y_n') = \text{cov}(\tilde{X}_{n,1})$, we have, $Q_n(Y_n')$ and $f(Y_n')w_n$ are independent and Gaussian. This implies

$$(8.3) \quad E(\|Y_n'\|^2) = E(\|\tilde{X}_{n,1}\|^2) = E(\|Q_n(\tilde{X}_{n,1})\|^2) + \|w_n\|^2\tilde{\sigma}_{f,n}^2,$$

hence, $\|w_n\|^2 \leq 4(M/\tilde{\sigma}_{f,n}^3)^{2/3} \leq 8(M/\delta^3)^{2/3}$ if $4M\delta^{-3}/n^{1/2} \leq 1$.

As $\|\tilde{S}_n/n^{1/2}\|^2 \leq \|Q_n(\tilde{S}_n/n^{1/2})\|^2 + 2\|Q_n(\tilde{S}_n/n^{1/2})\|\|f(\tilde{S}_n/n^{1/2})\|\|w_n\| + \|w_n\|^2|f(\tilde{S}_n/n^{1/2})|^2$, it follows that

$$
\begin{aligned}
E[\exp(-\lambda\rho_n f(\tilde{S}_n/n^{1/2}))&I\{\|\tilde{S}_n/n^{1/2}\|^2 \leq (2-5\varepsilon_n)\rho_n f(\tilde{S}_n/n^{1/2})\}] \\
&\geq E[\exp(-\lambda\rho_n f(\tilde{S}_n/n^{1/2})) \\
&\quad \times I\{\|Q_n(\tilde{S}_n/n^{1/2})\|^2 \leq (2-7\varepsilon_n)\rho_n f(\tilde{S}_n/n^{1/2})\}] \\
&\quad - P(\|w_n\|^2|f(\tilde{S}_n/n^{1/2})| > \varepsilon_n\rho_n) \\
&\quad - P(2\|w_n\|\|Q_n(\tilde{S}_n/n^{1/2})\| \geq \varepsilon_n\rho_n).
\end{aligned}
$$
(8.4)

Using Chebyshev's inequality along with (8.2) and (8.3), we have

$$(8.5) \quad \begin{aligned} P(\|w_n\|^2|f(\tilde{S}_n)/n^{1/2}| &> \varepsilon_n\rho_n) \\ &\leq \|w_n\|^4\tilde{\sigma}_{f,n}^2\varepsilon_n^{-2}\rho_n^{-2} \leq 4M^{2/3}\|w_n\|^2\varepsilon_n^{-2}\rho_n^{-2}, \end{aligned}$$

where $\|w_n\|$ is bounded as following (8.3). Likewise, it follows that

$$(8.6) \quad P(2\|w_n\|\|Q_n(\tilde{S}_n/n^{1/2})\| \geq \varepsilon_n\rho_n) \leq 16M^{2/3}\|w_n\|^2\varepsilon_n^{-2}\rho_n^{-2}.$$

Assuming that $\varepsilon_n\rho_n^{1/2} \to \infty$, we see that these two probabilities are of order $o(\rho_n^{-1})$.

(iii) Before we can proceed with the proof we need further lemmas.

LEMMA 13. *Let $Z_1$ and $Z_2$ be independent random variables and $c, d > 0$. Then*

$$(8.7) \quad P(Z_1 \geq cZ_2, Z_2 \geq 0) \leq r_2(d)[P(Z_1 \geq 0) + E(Z_1^+)/(cd)],$$

*where $r_2(d) = \sup_{x \geq 0} P(x \leq Z_2 < x + d)$.*

PROOF. Using the independence of $Z_1$ and $Z_2$, it follows that

$$P(Z_1 \geq cZ_2, Z_2 \geq 0)$$

$$= \sum_{j=1}^{\infty} P(Z_1 \geq cZ_2, (j-1)d \leq Z_2 < jd)$$

$$\leq \sum_{j=1}^{\infty} P(Z_1 \geq (j-1)cd)P((j-1)d \leq Z_2 < jd)$$

$$\leq r_2(d) \sum_{j=1}^{\infty} P(Z_1/(cd) \geq j-1) \leq r_2(d)[P(Z_1 \geq 0) + E(Z_1^+)/(cd)],$$

and the lemma is proved. □

LEMMA 14. *Let $V_n = \|Q_n(\tilde{S}_n/n^{1/2})\|^2 + \alpha_n G_1$, where $G_1$ is a standard normal random variable independent of $X_{n,1}, \ldots, X_{n,n}$, $\alpha_n \to 0$. If $\varepsilon_n' = 7\varepsilon_n + \alpha_n \log(1/\alpha_n)$, then*

$$E[\exp(-\lambda \rho_n f(\tilde{S}_n/n^{1/2})) I\{\|Q_n(\tilde{S}_n/n^{1/2})\|^2 \leq (2 - 7\varepsilon_n)\rho_n f(\tilde{S}_n/n^{1/2})\}]$$
$$\geq E[\exp(-\lambda \rho_n f(\tilde{S}_n/n^{1/2})) I\{V_n \leq (2 - \varepsilon_n')\rho_n f(\tilde{S}_n/n^{1/2}), f(\tilde{S}_n) > 0\}]$$
$$- o(\rho_n^{-1}).$$

PROOF. As

$$\{\|Q_n(\tilde{S}_n/n^{1/2})\|^2 \leq (2 - 7\varepsilon_n)\rho_n f(\tilde{S}_n/n^{1/2})\}$$
$$\supset \{V_n \leq (2 - \varepsilon_n')\rho_n f(\tilde{S}_n/n^{1/2})\} \cap \{f(\tilde{S}_n) > 0\}$$
$$\cap \{G_1 \geq -\log(1/\alpha_n)\rho_n f(\tilde{S}_n/n^{1/2})\},$$

and $G_1$ is symmetric, it is enough to show that

$$P(G_1 \geq \log(1/\alpha_n)\rho_n f(\tilde{S}_n/n^{1/2}), f(\tilde{S}_n) > 0) = o(\rho_n^{-1}).$$

Arguing as in part (i) (when estimating $p_n$), we see that

(8.8) $\quad r_n(d) = \sup_x P(x \leq f(\tilde{S}_n/n^{1/2}) \leq x + d) \leq d\delta^{-1} + C_{\mathrm{BE}} 8^{3/2} M \delta^{-3}/n^{1/2}.$

Applying Lemma 13 with $Z_1 = G_1, Z_2 = f(\tilde{S}_n/n^{1/2}), c_n = \log(1/\alpha_n)\rho_n, d_n = 1/c_n$, we find that the above probability is $\leq r_n(d_n)(1 + E[|G_1|]/2) = o(\rho_n^{-1})$. □

LEMMA 15. *Let $W_n = f(\tilde{S}_n/n^{1/2}) + \alpha_n' G_2$, where $G_2$ is a standard normal random variable independent of $X_{n,1}, \ldots, X_{n,n}, G_1$ and $\beta_n := \alpha_n' \rho_n \to 0$. If $\varepsilon_n'' = \varepsilon_n' + \beta_n^{1/2}$, we have*

$$E[\exp(-\lambda \rho_n f(\tilde{S}_n/n^{1/2})) I\{V_n \leq (2 - \varepsilon_n')\rho_n f(\tilde{S}_n/n^{1/2}), f(\tilde{S}_n) > 0\}]$$
$$\geq \exp(-\lambda^2 \beta_n^2/2)(E[\exp(-\lambda \rho_n W_n)$$
$$\times I\{V_n \leq (2 - \varepsilon_n'')\rho_n W_n, W_n > 0\}] - o(\rho_n^{-1})).$$

PROOF. By independence we obviously have

$$E[\exp(-\lambda\rho_n W_n)I\{V_n \le (2-\varepsilon'_n)\rho_n f(\tilde{S}_n/n^{1/2}), f(\tilde{S}_n) > 0\}]$$
$$= \exp(\lambda^2 \beta_n^2/2) E[\exp(-\lambda \rho_n f(\tilde{S}_n/n^{1/2}))$$
$$\times I\{V_n \le (2-\varepsilon'_n)\rho_n f(\tilde{S}_n/n^{1/2}), f(\tilde{S}_n) > 0\}].$$

We further have $\{V_n \le (2-\varepsilon'_n)\rho_n f(\tilde{S}_n/n^{1/2})\} =: A_n \supset B_n \cap C_n$, where

$$B_n = \{V_n \le (2-\varepsilon'_n - \beta_n^{1/2})\rho_n W_n\}, \qquad C_n = \{2\alpha'_n G_2 \le \beta_n^{1/2} f(\tilde{S}_n/n^{1/2})\}.$$

Therefore,

(8.9)
$$E[\exp(-\lambda\rho_n W_n)I_{A_n \cap \{f(\tilde{S}_n)>0\}}]$$
$$\ge E[\exp(-\lambda\rho_n W_n)I_{B_n \cap \{f(\tilde{S}_n)>0\}}] - P(C_n^c \cap \{f(\tilde{S}_n) > 0\}).$$

To bound the above probability we use once more Lemma 13. Setting $d_n = 1/c_n = \beta_n^{1/2}\rho_n^{-1}$, it follows that

(8.10) $\quad P(2\alpha'_n G_2 \ge \beta_n^{1/2} f(\tilde{S}_n/n^{1/2}), f(\tilde{S}_n) > 0) \le r_n(d_n)(1/2 + 2E[G_2^+]).$

Recalling (8.8), we see that $r_n(d_n) = o(\rho_n^{-1})$ so that it suffices to show that

(8.11)
$$E[\exp(-\lambda\rho_n W_n)I_{B_n \cap \{f(\tilde{S}_n)>0\}}]$$
$$\ge E[\exp(-\lambda\rho_n W_n)I_{B_n \cap \{W_n>0\}}] - o(\rho_n^{-1}).$$

To that end we first note that

$$E[\exp(-\lambda\rho_n W_n)I_{B_n \cap \{f(\tilde{S}_n)>0\}}]$$
$$\ge E[\exp(-\lambda\rho_n W_n)I_{B_n \cap \{W_n>0\} \cap \{f(\tilde{S}_n)>0\}}]$$
$$= E[\exp(-\lambda\rho_n W_n)I_{B_n \cap \{W_n>0\}}]$$
$$\quad - E[\exp(-\lambda\rho_n W_n)I_{B_n \cap \{W_n>0\} \cap \{f(\tilde{S}_n)\le 0\}}]$$
$$\ge E[\exp(-\lambda\rho_n W_n)I_{B_n \cap \{W_n>0\}}] - P(W_n > 0, f(\tilde{S}_n) \le 0).$$

Next, observe that

$$P(W_n > 0, f(\tilde{S}_n) \le 0) \le P(\alpha'_n G_2 \ge -f(\tilde{S}_n/n^{1/2}), -f(\tilde{S}_n) \ge 0),$$

which in view of Lemma 13 is bounded above by $r'_n(\alpha'_n)(1/2 + E[G_2^+])$, where $r'_n$ is defined as $r_n$ with $f$ replaced by $-f$. It is obvious that the upper bound in (8.8) also applies to $r'_n$ and we see that the above probability is of order $O(\alpha'_n) = o(\rho_n^{-1})$. This shows that (8.11) holds and Lemma 15 has been proven.

□

(iv) Recall that $Y'_n$ is a centered Gaussian random vector with covariance equal to that of $\tilde{X}_{n,1}$. Assuming that $Y'_n$ is independent of $G_1$, $G_2$, we set

$$\overline{V}_n = \|Q_n(Y'_n)\|^2 + \alpha_n G_1, \tag{8.12}$$
$$\overline{W}_n = f(Y'_n) + \beta_n \rho_n^{-1} G_2.$$

The purpose of this part of the proof is to show that

$$E[\exp(-\lambda \rho_n W_n) I\{V_n \leq (2 - \varepsilon''_n)\rho_n W_n, W_n > 0\}] \tag{8.13}$$
$$= E[\exp(-\lambda \rho_n \overline{W}_n) I\{\overline{V}_n \leq (2 - \varepsilon''_n)\rho_n \overline{W}_n, \overline{W}_n > 0\}] + o(\rho_n^{-1}).$$

To that end we first prove the following lemma.

LEMMA 16. *We have*

$$E[\exp(-\lambda \rho_n W_n) I\{W_n > 0\}] = E[\exp(-\lambda \rho_n \overline{W}_n) I\{\overline{W}_n > 0\}] + O(n^{-1/2}).$$

PROOF. Integration by parts yields that

$$E[(1 - \exp(-\lambda \rho_n W_n)) I\{W_n > 0\}] = \lambda \rho_n \int_0^\infty \exp(-\lambda \rho_n u) P\{W_n \geq u\}\, du.$$

Using the corresponding formula for $E[(1 - \exp(-\lambda \rho_n \overline{W}_n)) I\{\overline{W}_n > 0\}]$, we readily obtain that

$$|E[\exp(-\lambda \rho_n W_n) I\{W_n > 0\}] - E[\exp(-\lambda \rho_n \overline{W}_n) I\{\overline{W}_n > 0\}]|$$
$$\leq \lambda \rho_n \int_0^\infty \exp(-\lambda \rho_n u) |P\{W_n \geq u\} - P\{\overline{W}_n \geq u\}|\, du$$
$$+ |P\{W_n > 0\} - P\{\overline{W}_n > 0\}|,$$

which is obviously $\leq 2 \sup_u |P\{W_n \geq u\} - P\{\overline{W}_n \geq u\}|$. By conditioning on the independent variable $G_2$, we see that the last term in turn is bounded above by

$$\sup_x |P\{f(\tilde{S}_n/n^{1/2}) \geq x\} - P\{f(Y'_n) \geq x\}| \leq 2 C_{\text{BE}} 8^{3/2} M \delta^{-3}/n^{1/2},$$

where we have used once more the Berry–Esseen inequality. □

In view of Lemma 16 it is clear that (8.13) is proven once we have established the subsequent lemma.

LEMMA 17. *We have*

$$E[\exp(-\lambda \rho_n W_n) I\{0 < (2 - \varepsilon''_n)\rho_n W_n < V_n\}] \tag{8.14}$$
$$= E[\exp(-\lambda \rho_n \overline{W}_n) I\{0 < (2 - \varepsilon''_n)\rho_n \overline{W}_n < \overline{V}_n\}] + o(\rho_n^{-1}).$$

PROOF. We first note that

$$P(V_n \geq x) \leq P(G_1 \geq x/2) + P(\|Q_n(\tilde{S}_n)\| > (nx/2)^{1/2}),$$

and choosing $c > 0$ sufficiently large, we have from (8.2) and the Fuk–Nagaev type inequality presented in [11], page 338, that for $n$ large

$$P(V_n \geq c\log\rho_n)$$
$$\leq \exp\left\{-\frac{(c\log\rho_n)^2}{8}\right\} + \frac{72 \cdot 2^{11} nMq_n^3}{(cn\log\rho_n/2)^{3/2}} + \exp\left\{-\frac{cn\log\rho_n}{768nM^{2/3}}\right\}$$
$$= o(\rho_n^{-1}),$$

where we have used that $E(\|Q_n(\tilde{X}_{n,1})\|^2) \leq 4M^{2/3}$, which follows from (8.2) and (8.3). The latter relation also implies that $q_n = \|Q_n\| \leq 1 + \|w_n\|$ is bounded. Therefore,

(8.15)
$$E[\exp(-\lambda\rho_n W_n)I\{0 < (2 - \varepsilon_n'')\rho_n W_n < V_n\}]$$
$$= E[\exp(-\lambda\rho_n W_n)I\{(V_n, W_n) \in A_n\}] + o(\rho_n^{-1}),$$

where $A_n = \{(v,w) : 0 \leq w < (2 - \varepsilon_n'')^{-1} v/\rho_n \leq (2 - \varepsilon_n'')^{-1} c(\log\rho_n)/\rho_n\}$.

By an obvious modification of the above argument we find that also

(8.16)
$$E[\exp(-\lambda\rho_n \overline{W}_n)I\{0 < (2 - \varepsilon_n'')\rho_n \overline{W}_n < \overline{V}_n\}]$$
$$= E[\exp(-\lambda\rho_n \overline{W}_n)I\{(\overline{V}_n, \overline{W}_n) \in A_n\}] + o(\rho_n^{-1}).$$

It thus suffices to prove that

(8.17)
$$E[\exp(-\lambda\rho_n W_n)I\{(V_n, W_n) \in A_n\}]$$
$$= E[\exp(-\lambda\rho_n \overline{W}_n)I\{(\overline{V}_n, \overline{W}_n) \in A_n\}] + o(\rho_n^{-1}).$$

Let $f_{n,1}$ be the (two-dimensional) Lebesgue density function of $(V_n, W_n)$ and $f_{n,2}$ that of $(\overline{V}_n, \overline{W}_n)$. (These exist because we added an independent normal random vector.) Then we obviously have that the absolute difference of the two last expectations is bounded above by

$$\|f_{n,1} - f_{n,2}\|_\infty \operatorname{Area}(A_n).$$

Since $A_n$ is the triangle in the $(v,w)$ plane with base $c\log\rho_n$ and height $(2 - \varepsilon_n'')^{-1} c(\log\rho_n)/\rho_n$, we obviously have if $\varepsilon_n'' \leq 1$,

$$\operatorname{Area}(A_n) \leq c^2(\log\rho_n)^2/(2\rho_n),$$

and relation (8.17) immediately follows from the subsequent Lemma 18. □

LEMMA 18. *If $f_{n,1}$, $f_{n,2}$ are as above, where $\alpha_n = \beta_n = (\log\rho_n)^{-1}$, then we have for some $\gamma > 0$,*

(8.18) $$\|f_{n,1} - f_{n,2}\|_\infty = O(\rho_n^{-\gamma}).$$

PROOF. First, observe that by the inversion formula it is enough to show the characteristic functions $\phi_{n,1}$ and $\phi_{n,2}$ of $(V_n, W_n)$ and $(\overline{V}_n, \overline{W}_n)$, respectively, satisfy

$$\iint_{\mathbb{R}^2} |\phi_{n,1}(s,t) - \phi_{n,2}(s,t)|\, ds\, dt = O(\rho_n^{-\gamma}). \tag{8.19}$$

To verify (8.19) let

$$I_{n,1} = \int_{|s| \leq (\log \rho_n)^2} \int_{|t| \leq n^\tau} |\phi_{n,1}(s,t) - \phi_{n,2}(s,t)|\, ds\, dt,$$

$$I_{n,2} = \int_{|s| \leq (\log \rho_n)^2} \int_{n^\tau \leq |t| \leq \rho_n(\log \rho_n)^2} |\phi_{n,1}(s,t) - \phi_{n,2}(s,t)|\, ds\, dt,$$

$$I_{n,3} = \int_{|s| \geq (\log \rho_n)^2} \int_{t \in \mathbb{R}} |\phi_{n,1}(s,t) - \phi_{n,2}(s,t)|\, ds\, dt,$$

$$I_{n,4} = \int_{s \in \mathbb{R}} \int_{|t| \geq \rho_n(\log \rho_n)^2} |\phi_{n,1}(s,t) - \phi_{n,2}(s,t)|\, ds\, dt,$$

where $\tau > 0$ will be specified later.

It is obviously enough to show for $k = 1, 2, 3, 4$ and some $\gamma > 0$, that

$$I_{n,k} = O(\rho_n^{-\gamma}). \tag{8.20}$$

*Proof of* (8.20) *when* $k = 1$. Let $Y_{n,1}, \ldots, Y_{n,n}$ be i.i.d. copies of $Y_n'$ and for $-\infty < s, t < \infty$ and $x \in H$, define

$$F(x) = \exp\{is\|Q_n(x)\|^2 + itf(x)\}.$$

Then $F(\cdot)$ is Frechét differentiable in $x$, and Taylor's formula with integral remainder, see [7], page 70, implies

$$F(x+h) = F(x) + DF(x)(h) + \tfrac{1}{2}D^2F(x)(h,h) \\ + \tfrac{1}{2}E[(1-\tau)^2 D^3(x+\tau h)(h,h,h)], \tag{8.21}$$

where $\tau$ is uniform on $[0,1]$ and $D^k F(x)$ is the $k$th derivative of $F$ at $x$. Thus, by a standard argument we can conclude that

$$|E(F(\tilde{S}_n/n^{1/2}) - F(Y_n'))| \leq \sum_{k=1}^n J_k, \tag{8.22}$$

where

$$J_k = |E(F(W_k + n^{-1/2}\tilde{X}_{n,k})) - E(F(W_k + n^{-1/2}Y_{n,k}))| \tag{8.23}$$

and

$$W_k = (\tilde{X}_{n,1} + \cdots + \tilde{X}_{n,k-1} + Y_{n,k+1} + \cdots + Y_{n,n})/n^{1/2}, \tag{8.24}$$

for $k = 1, 2, \ldots, n$. Recall that $Y_{n,k}, 1 \leq k \leq n$ are independent Gaussian random vectors with the same distribution as $Y'_n$, which can be chosen independently of the random vectors $\tilde{X}_{n,k}, 1 \leq k \leq n$. Using (8.21), we expand the terms in $J_k$ with $x = W_k$ and $h = \tilde{X}_{n,k}/n^{1/2}$ or $h = Y_{n,k}/n^{1/2}$, respectively. Since $\tilde{X}_{n,k}$ and $Y_{n,k}$ are independent of $W_k$ and $\tilde{X}_{n,k}$ and $Y_{n,k}$ both have mean zero with common covariance functions, the terms containing derivatives up to second order coincide so that

(8.25) $$J_k \leq n^{-3/2}(J'_k + J''_k),$$

where

$$J'_k = \tfrac{1}{2}|E((1-\tau)^2 D^3 F(W_k + \tau \tilde{X}_{n,k}/n^{1/2})(\tilde{X}_{n,k}, \tilde{X}_{n,k}, \tilde{X}_{n,k}))|$$

and

$$J''_k = \tfrac{1}{2}|E((1-\tau)^2 D^3 F(W_k + \tau Y_{n,k}/n^{1/2})(Y_{n,k}, Y_{n,k}, Y_{n,k}))|.$$

Since $\|f\| = 1$, we have

$$|DF^3(x)(h,h,h)|$$
$$\leq 12s^2 \|Q_n(h)\|^3 \|Q_n(x)\| + 6|st|\|Q_n(h)\|^2 \|h\|$$
$$+ 2^2(8|s|^3 \|Q_n(x)\|^3 \|Q_n(h)\|^3 + |t|^3 \|h\|^3),$$

and, therefore, since $0 \leq \tau \leq 1$ (and setting $q_n = \|Q_n\|$), we have

(8.26)
$$\begin{aligned}
J'_k &\leq \tfrac{1}{2}E[32|s|^3 \|Q_n(\tilde{X}_{n,k})\|^3 \|Q_n(W_k + \tau \tilde{X}_{n,k}/n^{1/2})\|^3 \\
&\quad + 6|st|\|Q_n(\tilde{X}_{n,k})\|^2 \|\tilde{X}_{n,k}\| + 4|t|^3 \|\tilde{X}_{n,k}\|^3 \\
&\quad + 12s^2 \|Q_n(\tilde{X}_{n,k})\|^3 \|Q_n(W_k + \tau \tilde{X}_{n,k}/n^{1/2})\|] \\
&\leq \tfrac{1}{2}E[32|s|^3 q_n^3 \|\tilde{X}_{n,k}\|^3 \{4(q_n^3 \|W_k\|^3 + q_n^3 \|\tilde{X}_{n,k}\|^3/n^{3/2})\} \\
&\quad + 6|st|q_n^2 \|\tilde{X}_{n,k}\|^3 + 4|t|^3 \|\tilde{X}_{n,k}\|^3 \\
&\quad + 12s^2 q_n^4 \|\tilde{X}_{n,k}\|^3 (\|W_k\| + \|\tilde{X}_{n,k}\|/n^{1/2})].
\end{aligned}$$

Now $P(\|\tilde{X}_{n,k}\| \leq 2\delta n^{1/2}) = 1$, $E(\|\tilde{X}_{n,k}\|^3) \leq 8M$, and $\tilde{X}_{n,k}$ and $W_k$ are independent, so

(8.27)
$$\begin{aligned}
J'_k &\leq [512 M q_n^6 |s|^3 E(\|W_k\|^3 + 8\delta^3) \\
&\quad + (24 M q_n^2 |st| + 16M|t|^3) + 48 M s^2 q_n^4 E(\|W_k\| + 2\delta)].
\end{aligned}$$

We need an upper bound for $E(\|W_k\|^3)$. To that end we first note that by convexity and that the $Y_{n,k}$'s are independent and identically distributed Gaussian random vectors, we have

(8.28) $$E\|W_k\|^3 \leq 4E(\|\tilde{S}_n\|^3)/n^{3/2} + 4E(\|Y_{n,1}\|^3).$$

Applying Proposition 6.8 of [18], we have

(8.29) $$E(\|\tilde{S}_n\|^3) \leq 2 \cdot 4^3 8 \delta^3 n^{3/2} + 2(4b_0)^3,$$

where

$$b_0 = \inf\left\{b : P\left(\max_{1 \leq j \leq n}\left\|\sum_{m=1}^j \tilde{X}_{n,m}\right\| > b\right) \leq (2 \cdot 4^3)^{-1}\right\}.$$

Using Proposition 1.1.2 in [10], we further have

$$P\left(\max_{1 \leq j \leq n}\left\|\sum_{m=1}^j \tilde{X}_{n,m}\right\| > b\right) \leq 3 \max_{1 \leq j \leq n} P\left(\left\|\sum_{m=1}^j \tilde{X}_{n,m}\right\| > b/3\right).$$

Since we are in a Hilbert space, we have

(8.30) $$\max_{1 \leq j \leq n} E\left(\left\|\sum_{m=1}^j \tilde{X}_{n,m}\right\|^2\right) \leq nE(\|\tilde{X}_{n,1}\|^2) \leq 4M^{2/3}n$$

and via Markov's inequality, it follows that $b_0 \leq 144 M^{1/3} n^{1/2}$.

Employing the trivial inequality $\delta^3 \leq M$, we see that

(8.31) $$E(\|\tilde{S}_n\|^3) \leq AMn^{3/2},$$

where $A$ is a universal constant.

Using the equivalence of the moments of Gaussian random variables (see, e.g., Corollary 3.2 of [18]), we also have

(8.32) $$E(\|Y_{n,1}\|^3) \leq A'E(\|\tilde{X}_{n,1}\|^2)^{3/2} \leq 8A'M,$$

where $A'$ is another universal constant. Since $E(\|W_k\|) \leq (E(\|W_k\|^3))^{1/3}$, by combining (8.28), (8.31) and (8.32) we have if $4M\delta^{-3}/n^{1/2} \leq 1$ and, consequently, $q_n \leq 1 + \|w_n\| \leq 1 + \sqrt{8}M^{1/3}\delta^{-1}$,

(8.33) $$J'_k \leq C_1[|s|^3 + |t|^3 + 1],$$

where $C_1$ is a finite constant depending only on $M$ and $\delta$.

Here we have used that $|st| \leq (s^2 + t^2)/2$, $s^2 \leq |s|^3 + 1$ and $t^2 \leq |t|^3 + 1$.

Similarly, it follows that if $4M\delta^{-3}/n^{1/2} \leq 1$, we have

(8.34) $$J''_k \leq C_2[|s|^3 + |t|^3 + 1],$$

where $C_2$ is another finite constant depending only on $M, \delta$.

Combining (8.22), (8.25), (8.33) and (8.34) with $C_3 = C_1 + C_2$, we see that under the assumption $4M\delta^{-3}/n^{1/2} \leq 1$,

(8.35) $$|E(F(\tilde{S}_n/n^{1/2}) - F(Y'_n))| \leq C_3 n^{-1/2}[|s|^3 + |t|^3 + 1].$$

Enlarging the constant $C_3$ if necessary, we finally see that this is also the case if $4M\delta^{-3}/n^{1/2} > 1$. [Use the trivial fact that $|E(F(\tilde{S}_n/n^{1/2}) - F(Y'_n))| \leq 2$.]

By independence we obviously have

$$
\begin{aligned}
|\phi_{n,1}(s,t) - \phi_{n,2}(s,t)| \\
(8.36) \qquad = |E(e^{is\alpha_n G_1 + it\alpha'_n G_2})| \cdot |E(F(\tilde{S}_n/n^{1/2}) - F(Y'_n))| \\
\leq |E(F(\tilde{S}_n/n^{1/2}) - F(Y'_n))|.
\end{aligned}
$$

Thus, (8.35) and (8.36) imply

$$I_{n,1} \leq 4C_3 n^{-1/2} \int_0^{(\log \rho_n)^2} \left( \int_0^{n^\tau} (s^3 + t^3 + 1)\, dt \right) ds = O(n^{-1/2 + 4\tau}(\log n)^2),$$

which is of order $O(n^{-1/4}) = O(\rho_n^{-1/2})$, provided $0 < \tau < 1/16$. Thus, (8.20) holds for $k = 1$ with $\gamma = 1/2$ if $0 < \tau < 1/16$.

*Proof of (8.20) when $k = 2$.* Let $\Delta_n = \sum_{j=1}^{k_n} \tilde{X}_{n,j}/n^{1/2}$ and $U_n = \tilde{S}_n/n^{1/2} - \Delta_n$, where $k_n \leq n$ will be specified later. Then

$$\phi_{n,1}(s,t) = E(\exp[is\|Q_n(\Delta_n + U_n)\|^2 + itf(\Delta_n + U_n) + is\alpha_n G_1 + it\alpha'_n G_2]),$$

and we also define

$$\begin{aligned}
\widetilde{\phi}_{n,1}(s,t) = E(\exp[is\|Q_n(\Delta_n)\|^2 + 2is\langle Q_n(\Delta_n), \\
Q_n(U_n)\rangle + itf(\Delta_n + U_n) + is\alpha_n G_1 + it\alpha'_n G_2]).
\end{aligned}$$

Then, since $|e^{ix} - 1| \leq |x|$, we have because of (8.3),

$$
\begin{aligned}
(8.37) \quad |\phi_{n,1}(s,t) - \widetilde{\phi}_{n,1}(s,t)| &\leq E(|e^{is\|Q_n(U_n)\|^2} - 1|) \\
&\leq |s| E(\|Q_n(U_n)\|^2) \leq |s| E(\|\tilde{X}_{n,1}\|^2)\left(1 - \frac{k_n}{n}\right).
\end{aligned}
$$

Next observe that

$$
\begin{aligned}
(8.38) \quad |\widetilde{\phi}_{n,1}(s,t)|^2 &= |E(E((\cdots)|\Delta_n)|^2) \exp(-\alpha_n{}^2 s^2 - \alpha'_n{}^2 t^2) \\
&\leq E(|E((\cdots)|\Delta_n)|^2) \\
&= E(|E(\exp\{2is\langle Q_n(\Delta_n), Q_n(U_n)\rangle + itf(U_n)\}|\Delta_n)|^2),
\end{aligned}
$$

where $(\cdots) = \exp(is\|Q_n(\Delta_n)\|^2 + 2is\langle Q_n(\Delta_n), Q_n(U_n)\rangle + itf(\Delta_n + U_n))$. Thus, if $U'_n$ is an independent copy of $U_n$ which is also independent of $\Delta_n$, we readily obtain that

$$(8.39) \qquad |\widetilde{\phi}_{n,1}(s,t)|^2 \leq E(e^{2is\langle Q_n(\Delta_n), Q_n(U_n^*)\rangle + itf(U_n^*)}),$$

where $U_n^* = U_n - U'_n$ is the symmetrization of $U_n$.

Denoting the distribution of $Q_n(\Delta_n)$ by $\mu_n$, we have

$$E(\exp\{2is\langle Q_n(\Delta_n), Q_n(U_n^*)\rangle + itf(U_n^*)\})$$

(8.40)
$$= \int_H E[\exp\{2is\langle Q_n(U_n^*), x\rangle + itf(U_n^*)\}] \, d\mu_n(x)$$

$$= \int_H [E(\cos[(2s/n^{1/2})\langle Q_n(\tilde{X}_n^*), x\rangle + (t/n^{1/2})f(\tilde{X}_n^*)])]^{n-k_n} \, d\mu_n(x),$$

where $\tilde{X}_n^* = \tilde{X}_{n,1} - \tilde{X}_{n,2}$ is the symmetrization of $\tilde{X}_{n,1}$.

Now for each $x \in H$,

$$E(\cos[(2s/n^{1/2})\langle Q_n(\tilde{X}_n^*), x\rangle + (t/n^{1/2})f(\tilde{X}_n^*)])$$

$$\leq 1 - \tfrac{1}{2}E(((2s/n^{1/2})\langle Q_n(\tilde{X}_n^*), x\rangle + (t/n^{1/2})f(\tilde{X}_n^*))^2)$$

$$+ \tfrac{1}{6}E(|(2s/n^{1/2})\langle Q_n(\tilde{X}_n^*), x\rangle + (t/n^{1/2})f(\tilde{X}_n^*)|^3)$$

$$\leq 1 - (t^2/(2n))E(f^2(\tilde{X}_n^*)) + \tfrac{2}{3}(|t|^3/n^{3/2})E(|f(\tilde{X}_n^*)|^3)$$

$$+ \tfrac{16}{3}(|s|^3/n^{3/2})E(\|Q_n(\tilde{X}_n^*)\|^3)\|x\|^3$$

$$\leq 1 - (t^2/n)E(f^2(\tilde{X}_{n,1})) + \tfrac{64}{3}M(|t|^3 + 8q_n^3\|x\|^3|s|^3)/n^{3/2},$$

where we first have used that $\langle Q_n(\tilde{X}_n^*), x\rangle$ and $f(\tilde{X}_n^*)$ are uncorrelated by our choice of $w_n$ and then relation (8.2), along with the convexity of the function $u \to |u|^3$.

If $4M\delta^{-3}/n^{1/2} \leq 1$ so that $E(f^2(\tilde{X}_{n,1})) \geq \delta^2/2$, we have for $\|x\| \leq |t|/(2|s|q_n)$,

$$E(\cos[(2s/n^{1/2})\langle Q_n(\tilde{X}_n^*), x\rangle + (t/n^{1/2})f(\tilde{X}_n^*)])$$

(8.41)
$$\leq 1 - t^2\delta^2/(2n) + \tfrac{128}{3}(|t|^3/n^{3/2})M$$

$$\leq 1 - t^2\delta^2/(4n),$$

provided $|t| \leq (3/512)\delta^2 n^{1/2}/M$.

Combining (8.37) with (8.39)–(8.41) and that $1 - x \leq e^{-x}$, we see that

$$|\phi_{n,1}(s,t)| \leq |\widetilde{\phi}_{n,1}(s,t)| + |s|E(\|\tilde{X}_{n,1}\|^2)(1 - k_n/n)$$

$$\leq |s|E(\|\tilde{X}_{n,1}\|^2)(1 - k_n/n)$$

(8.42)
$$+ (E(\exp\{2is\langle Q_n(\Delta_n), Q(U_n^*)\rangle + itf(U_n^*)\}))^{1/2}$$

$$\leq |s|E(\|\tilde{X}_{n,1}\|^2)(1 - k_n/n) + \exp\{-t^2\delta^2(n-k_n)/(8n)\}$$

$$+ (P\{\|Q_n(\Delta_n)\| \geq (2q_n)^{-1}|t/s|\})^{1/2},$$

if $|t| \leq (3/512)\delta^2 n^{1/2}/M$ and $|s| \leq (\log \rho_n)^2$.

Taking $k_n = n - [n/|t|^{3/2}] - 1$, we thus have by (8.3) and the Fuk–Nagaev inequality as given in [11] and such $s$, $t$ that

$$
\begin{aligned}
|\phi_{n,1}(s,t)| &\leq 4|s|M^{2/3}(|t|^{-3/2} + n^{-1}) + \exp\{-t^2\delta^2|t|^{-3/2}/8\} \\
&\quad + [9 \cdot 2^{17} M(|t/s|)^{-3} q_n^3]^{1/2} n^{-1/4} \\
&\quad + \exp\{-|t/s|^2/(768 E(\|\tilde{X}_{n,1}\|^2))\} \\
&\leq 4|s|M^{2/3}(|t|^{-3/2} + n^{-1}) + \exp\{-|t|^{1/2}\delta^2/8\} \\
&\quad + C_4\{|s/t|^{3/2} n^{-1/4} + \exp\{-C_5|t/s|^2\}\}.
\end{aligned}
$$
(8.43)

In the first inequality of (8.43) when we apply the Fuk–Nagaev inequality, we use the fact that for $I_{n,2}$, the ratio $|t/s| \geq n^\tau/(\log \rho_n)^2$, and that $E(\|\Delta_n\|) \leq E(\|\tilde{S}_n\|^2)^{1/2}/n^{1/2} \leq M^{1/3} < \infty$. Also, recall that $q_n \leq 1 + \sqrt{8}M^{1/3}/\delta$ if $n$ is large enough.

As $f(Y_n'), \|Q_n(Y_n')\|, G_1$ and $G_2$ are independent random variables, we readily obtain for $-\infty < s, t < \infty$ (assuming $4M\delta^{-3}/n^{1/2} \leq 1$),

$$
\begin{aligned}
|\phi_{n,2}(s,t)| &\leq E(\exp(itf(Y_n'))) = \exp(-t^2 E(f^2(\tilde{X}_{n,1})/2)) \\
&\leq \exp(-t^2\delta^2/4),
\end{aligned}
$$
(8.44)

which is for $|t| \geq 1$ dominated by $\exp\{-|t|^{1/2}\delta^2/8\}$. It thus follows that for large $n$,

$$
\begin{aligned}
|\phi_{n,1}(s,t) &- \phi_{n,2}(s,t)| \\
&\leq C_6(|s||t|^{-3/2} + |s|n^{-1}) + C_7 \exp\{-C_8|t|^{1/2}\} \\
&\quad + C_4|s/t|^{3/2} n^{-1/4} + C_4 \exp\{-C_5|t/s|^2\},
\end{aligned}
$$
(8.45)

provided that $|s| \leq (\log \rho_n)^2$ and $n^\tau \leq |t| \leq (3/512)\delta^2 n^{1/2}/M$.

The constants $C_i$ depend on $M$ and $\delta$ only and are strictly positive and finite.

Integrating over the region related to $I_{n,2}$, we have $I_{n,2} = O(n^{-\tau/2}(\log n)^4)$. Thus, (8.20) holds for $k = 2$ with $\gamma = 0.9\tau$ (say).

*Proof of* (8.20) *when* $k = 3, 4$. Recalling the definition of $V_n$, $W_n$, $\overline{V}_n$ and $\overline{W}_n$ and that $\alpha_n = (\log \rho_n)^{-1}$, $\alpha_n' = (\rho_n \log \rho_n)^{-1}$, we see that

$$
\begin{aligned}
I_{n,3} &\leq 2 \int_{|s| \geq (\log \rho_n)^2} \int_{t \in \mathbb{R}} \exp(-\alpha_n^2 s^2/2 - \alpha_n'^2 t^2/2)\, ds\, dt \\
&\leq \frac{4}{\log \rho_n} \exp(-(\log \rho_n)^2/2)\sqrt{2\pi} \rho_n \log \rho_n,
\end{aligned}
$$

which is obviously of order $o(\rho_n^{-1})$. A similar calculation shows finally that $I_{n,4} = o(\rho_n^{-1})$, thereby completing the proof of Lemma 18. $\square$

(vi) Given (8.13) we now investigate the asymptotic behavior of
$$E[\exp(-\lambda \rho_n \overline{W}_n) I\{\overline{V}_n \leq (2 - \varepsilon_n'') \rho_n \overline{W}_n, \overline{W}_n > 0\}].$$
We first show that we can remove the smoothing variable $\alpha_n G_1$. Arguing as in the proof of Lemma 14, we find that
$$E[\exp(-\lambda \rho_n \overline{W}_n) I\{\overline{V}_n \leq (2 - \varepsilon_n'') \rho_n \overline{W}_n, \overline{W}_n > 0\}]$$
$$\geq E[\exp(-\lambda \rho_n \overline{W}_n) I\{\|Q_n(Y_n')\|^2 \leq (2 - 2\varepsilon_n'') \rho_n \overline{W}_n\}]$$
$$- P\{\alpha_n G_1 \geq \varepsilon_n'' \rho_n \overline{W}_n, \overline{W}_n > 0\},$$
where we have $P\{\alpha_n G_1 \geq \varepsilon_n'' \rho_n \overline{W}_n, \overline{W}_n > 0\} = o(\rho_n^{-1})$ by Lemma 13. [Recall that $\varepsilon_n''/\alpha_n \to \infty$ and use the fact that the densities of the random variables $\overline{W}_n \sim \text{normal}(0, \tilde{\sigma}_{f,n}^2 + \alpha_n^2)$ are uniformly bounded.]

Let $g_{n,1}$ and $g_{n,2}$ be the (normal) densities of $f(Y_n')$ and $\overline{W}_n$, respectively. Then, using the inversion formula for densities, we see that
$$(8.46) \quad \sqrt{2\pi}\|g_{n,1} - g_{n,2}\|_\infty = \tilde{\sigma}_{f,n}^{-1} - (\tilde{\sigma}_{f,n}^2 + {\alpha_n'}^2)^{-1/2} = o(\rho_n^{-2}).$$
Let further $\nu_n$ be the distribution of $\|Q_n(Y_n')\|^2$. By independence of the variables $\|Q_n(Y_n')\|^2, f(Y_n')$, and $G_2$, we have then
$$(8.47) \quad \begin{aligned} &E[\exp(-\lambda \rho_n \overline{W}_n) I\{\|Q_n(Y_n')\|^2 \leq (2 - 2\varepsilon_n'') \rho_n \overline{W}_n\}] \\ &= \int_0^\infty \int_{x/[(2-\varepsilon_n'')\rho_n]}^\infty \exp(-\lambda \rho_n z) g_{n,2}(z) \, dz \, d\nu_n(x) \end{aligned}$$
and
$$(8.48) \quad \begin{aligned} &E[\exp(-\lambda \rho_n f(Y_n')) I\{\|Q_n(Y_n')\|^2 \leq (2 - 2\varepsilon_n'') \rho_n f(Y_n')\}] \\ &= \int_0^\infty \int_{x/[(2-\varepsilon_n'')\rho_n]}^\infty \exp(-\lambda \rho_n z) g_{n,1}(z) \, dz \, d\nu_n(x). \end{aligned}$$
Combining (8.46)–(8.48), we can infer that
$$(8.49) \quad \begin{aligned} &E[\exp(-\lambda \rho_n \overline{W}_n) I\{\|Q_n(Y_n')\|^2 \leq (2 - 2\varepsilon_n'') \rho_n \overline{W}_n\}] \\ &= E[\exp(-\lambda \rho_n f(Y_n')) I\{\|Q_n(Y_n')\|^2 \leq (2 - 2\varepsilon_n'') \rho_n f(Y_n')\}] + o(\rho_n^{-2}). \end{aligned}$$

(vii) Next, observe that
$$E[\exp(-\lambda \rho_n f(Y_n')) I\{\|Q_n(Y_n')\|^2 \leq 2\rho_n f(Y_n')\}]$$
$$- E[\exp(-\lambda \rho_n f(Y_n')) I\{\|Q_n(Y_n')\|^2 \leq (2 - 2\varepsilon_n'') \rho_n f(Y_n')\}]$$
$$= \int_0^\infty \int_{x/(2\rho_n)}^{x/[(2-2\varepsilon_n'')\rho_n]} \exp(-\lambda \rho_n z) g_{n,1}(z) \, dz \, d\nu_n(x)$$
$$\leq \frac{\varepsilon_n''}{2(1 - \varepsilon_n'')} \|g_{n,1}\|_\infty E[\|Q_n(Y_n')\|]/\rho_n$$
$$= o(\rho_n^{-1}).$$

By independence we have, for any $A > 0$,

$$E[\exp(-\lambda \rho_n f(Y'_n))I\{\|Q_n(Y'_n)\|^2 \leq 2\rho_n f(Y'_n)\}]$$
$$\geq \exp(-2\lambda A)P\{A < f(Y'_n)\rho_n < 2A\}P\{\|Q_n(Y'_n)\|^2 \leq A\},$$

which in turn via Markov's inequality and (8.2) and (8.3) is greater than or equal to

$$\exp(-2\lambda A)P\{A < f(Y'_n)\rho_n < 2A\}/2,$$

if we choose $A = 8M^{2/3}$. The density functions of $f(Y'_n)$ are eventually uniformly positive in a neigborhood of zero so that

(8.50) $\quad \liminf_{n \to \infty} \rho_n E[\exp(-\lambda \rho_n f(Y'_n))I\{\|Q_n(Y'_n)\|^2 \leq 2\rho_n f(Y'_n)\}] > 0,$

and we can conclude that as $n \to \infty$,

(8.51) $\quad \liminf_{n \to \infty} I_n/E[\exp(-\lambda \rho_n f(Y'_n))I\{\|Q_n(Y'_n)\|^2 \leq 2\rho_n f(Y'_n)\}] \geq 1.$

(viii) Given $n \geq 1$, let $Y''_n$ be a Gaussian mean zero random vector which is independent of $Y'_n$ so that

$$\mathcal{L}(Y_n) = \mathcal{L}(Y'_n + Y''_n), \qquad n \geq 1.$$

[Such a sequence exists since $\operatorname{cov}(Y_n) - \operatorname{cov}(Y'_n)$ is positive semidefinite, as can easily be seen from the definition of these random vectors.]

Denoting the density function of $f(Y_n)$ by $g_n$, it follows that $\|g_{n,1} - g_n\|_\infty \to 0$, which in turn by the independence of $f(Y_n)$ and $Q_n(Y'_n)$ and a slight modification of (8.48) implies

(8.52) 
$$E[\exp(-\lambda \rho_n f(Y'_n))I\{\|Q_n(Y'_n)\|^2 \leq 2\rho_n f(Y'_n)\}]$$
$$= E[\exp(-\lambda \rho_n f(Y_n))I\{\|Q_n(Y'_n)\|^2 \leq 2\rho_n f(Y_n)\}] + o(\rho_n^{-1}).$$

Setting

$$v_n = E[f(Y''_n)Y''_n]/E[f^2(Y''_n)]^{1/2} \quad \text{and} \quad Q''_n(x) = x - v_n f(x), \qquad x \in H,$$

if $E[f^2(Y''_n)] > 0$, we obviously have

$$\|Y_n\|^2 = \|Q_n(Y'_n) + Q'_n(Y''_n) + f(Y'_n)w_n + f(Y''_n)v_n\|^2,$$

where the variables $Q_n(Y'_n), Q'_n(Y''_n), f(Y'_n)$ and $f(Y''_n)$ are independent. It thus follows that

(8.53)
$$E[\exp(-\lambda \rho_n f(Y_n))I\{\|Y_n\|^2 \leq 2\rho_n f(Y_n)\}]$$
$$= \int_{-\infty}^{\infty} \int_{-z_1}^{\infty} e^{-\lambda(z_1+z_2)} p_n(z_1, z_2) g_{n,3}(z_2) \, dz_2 \, g_{n,1}(z_1) \, dz_1,$$

where $p_n(z_1, z_2) = P\{\|Q_n(Y'_n) + Q'_n(Y''_n) + z_1 w_n + z_2 v_n\|^2 \leq 2\rho_n(z_1 + z_2)\}$ and $g_{n,3}$ is the density of $f(Y''_n)$.

By the inequality of Anderson we have, for $z_1, z_2 \in \mathbb{R}$,

$$p_n(z_1, z_2) \leq P\{\|Q_n(Y_n') + Q_n'(Y_n'')\|^2 \leq 2\rho_n(z_1 + z_2)\}$$
$$\leq P\{\|Q_n(Y_n')\|^2 \leq 2\rho_n(z_1 + z_2)\},$$

which in combination with (8.53) implies

$$E[\exp(-\lambda \rho_n f(Y_n))I\{\|Y_n\|^2 \leq 2\rho_n f(Y_n)\}]$$
$$\leq E[\exp(-\lambda \rho_n f(Y_n))I\{\|Q_n(Y_n')\|^2 \leq 2\rho_n f(Y_n)\}].$$

Recalling (8.50) and (8.52), we obtain the desired result. $\square$

**Acknowledgments.** We would like to thank Professor Alex de Acosta, whose remarks led us to begin our investigation of these problems and Professor Vidas Bentkus for informing us about very recent work on Berry–Esseen type results for $U$-statistics. Thanks are also due to an anonymous referee for providing a reference for the proof of Lemma 11. An alternative proof had been given to us by Professor Pieter de Groen at an earlier stage which was appreciated as well.

DEPARTEMENT WISKUNDE
VRIJE UNIVERSITEIT BRUSSEL
PLEINLAAN 2
B-1050 BRUSSEL
BELGIUM
E-MAIL: ueinmahl@vub.ac.be

DEPARTMENT OF MATHEMATICS
UNIVERSITY OF WISCONSIN
MADISON, WISCONSIN 53706
USA
E-MAIL: kuelbs@math.wisc.edu